\documentclass[10pt,a4paper]{article}
\usepackage{amsmath,amssymb,amsthm,graphics,graphicx,epic,eepic,multicol,ascmac}
\usepackage[all]{xy}

\numberwithin{equation}{section}
\theoremstyle{theorem}
\newtheorem{thm}{Theorem}[section]
\newtheorem{prop}[thm]{Proposition}
\newtheorem{lem}[thm]{Lemma}
\newtheorem{rem}[thm]{Remark}

\newtheorem{ex}[thm]{Example}

\theoremstyle{definition}
\newtheorem{defn}[thm]{Definition}

\def\al{\alpha}

\def\ep{\epsilon}
\def\wht(#1){\widehat{\ #1\ }}

\newcommand{\cA}{{\mathcal A}}

\newcommand{\cF}{{\mathcal F}}

\newcommand{\frg}{\mathfrak g}
\newcommand{\frh}{\mathfrak h}

\newcommand{\frn}{\mathfrak n}

\newcommand{\bbC}{\mathbb C}

\newcommand{\ch}{\mathrm{ch}}

\newcommand{\lbr}{\begin{bmatrix}}
\newcommand{\rbr}{\end{bmatrix}}

\newcommand{\cd}{commutative diagram }

\def\ge{\frg}

\def\cP{{\mathcal P}}
\def\al{\alpha}

\def\beneme{\begin{enumerate}}
\def\beq{\begin{equation}}
\def\beqn{\begin{eqnarray}}
\def\beqnn{\begin{eqnarray*}}
\def\bfi{{\mathbf i}}
\def\bfii0{{\bf i_0}}

\def\bbra#1,#2,#3{\left\{\begin{array}{c}\hspace{-5pt}
#1;#2\\ \hspace{-5pt}#3\end{array}\hspace{-5pt}\right\}}
\def\cd{\cdots}
\def\ci(#1,#2){c_{#1}^{(#2)}}
\def\Ci(#1,#2){C_{#1}^{(#2)}}
\def\mpp(#1,#2,#3){#1^{(#2)}_{#3}}
\def\bCi(#1,#2){\ovl C_{#1}^{(#2)}}
\def\ch(#1,#2){c_{#2,#1}^{-h_{#1}}}
\def\cc(#1,#2){c_{#2,#1}}

\def\del{\delta}
\def\Del{\Delta}

\def\di(#1,#2){D_{#1}^{(#2)}}
\def\dbi(#1,#2){\ovl D_{#1}^{(#2)}}

\def\eneme{\end{enumerate}}
\def\ep{\epsilon}
\def\eeq{\end{equation}}
\def\eeqn{\end{eqnarray}}
\def\eeqnn{\end{eqnarray*}}

\def\gau#1,#2{\left[\begin{array}{c}\hspace{-5pt}#1\\
\hspace{-5pt}#2\end{array}\hspace{-5pt}\right]}

\def\ji(#1,#2){j_{#1}^{(#2)}}

\def\lan{\langle}

\def\lm{\lambda}
\def\Lm{\Lambda}

\def\nd{\noindent}

\def\ovl{\overline}

\def\qq{\qquad}
\def\q{\quad}
\def\qed{\hfill\framebox[2mm]{}}

\def\ran{\rangle}

\def\TY(#1,#2,#3){#1^{(#2)}_{#3}}

\def\UU{{\mathcal U}}

\def\xxi(#1,#2,#3){\displaystyle {}^{#1}\Xi^{(#2)}_{#3}}
\def\xsi(#1,#2,#3){\displaystyle {}^{#1}\Sigma^{(#2)}_{#3}}
\def\xE(#1,#2,#3){\displaystyle {}^{#1}E_{#2}[#3]}
\def\xF(#1,#2){\displaystyle {}^{#1}F_{#2}}
\def\xx(#1,#2){\displaystyle {}^{#1}\Xi_{#2}}
\def\W1{W(\varpi_1)}

\def\m@th{\mathsurround=0pt}
\def\fsquare(#1,#2){
\hbox{\vrule$\hskip-0.4pt\vcenter to #1{\normalbaselines\m@th
\hrule\vfil\hbox to #1{\hfill$\scriptstyle #2$\hfill}\vfil\hrule}$\hskip-0.4pt
\vrule}}

\newcommand{\ba}{\begin{array}}
\newcommand{\ea}{\end{array}}

\newcommand{\eq}{\begin{eqnarray}}
\newcommand{\eneq}{\end{eqnarray}}

\title{\textbf{\large{Explicit Forms of Cluster Variables 
on Double Bruhat Cells $G^{u,e}$ of type B}}}
\author{\normalsize{YUKI KANAKUBO\thanks{Division of Mathematics, 
Sophia University, Kioicho 7-1, Chiyoda-ku, Tokyo 102-8554,
Japan: {j\_chi\_sen\_you\_ky@eagle.sophia.ac.jp}
}}}
\date{}

\begin{document}

\maketitle
\vspace{-10pt}

\begin{abstract}
Let $G$
be a simply connected simple algebraic group over $\mathbb{C}$ of type ${\rm B}_r$, 
$B$ and $B_-$ be its two opposite Borel subgroups, and $W$ be the associated 
 Weyl group. 
For $u$, $v\in W$, 
it is known that the coordinate ring 
${\mathbb C}[G^{u,v}]$ of the double Bruhat cell $G^{u,v}=BuB\cap B_-vB_-$ is 
isomorphic to an upper cluster algebra $\ovl{\cA}(\textbf{i})_{{\mathbb C}}$ 
and generalized minors $\Delta(k;\textbf{i})$ are 
the cluster variables of ${\mathbb C}[G^{u,v}]$\cite{A-F-Z}. It is also shown that ${\mathbb C}[G^{u,v}]$ have a structure of cluster algebra \cite{GY}.
In the case $v=e$, 
we shall describe the generalized minor $\Delta(k;\textbf{i})$
 explicitly.
\end{abstract}


\tableofcontents

\section{Introduction}

Let $G$ be a simply connected simple algebraic group over $\bbC$ of rank
$r$, $B,B_-\subset G$ the opposite Borel subgroups, $H:=B\cap B_-$ the maximal torus, 
$N\subset B$, $N_-\subset B_-$ the maximal unipotent 
subgroups and $W$ the associated Weyl group. 

Fomin and Zelevinsky have invented cluster algebras for the study of total positivity and dual semi canonical bases in 2002 \cite{FZ2}. They are commutative algebras generated by so-called {\it cluster variables}. Choosing a part of the cluster variables properly, we can combinatorially calculate other variables from them (Sect.{\ref{CluSect}}). These chosen variables are called initial cluster variables.

In \cite{A-F-Z}, Berenstein, Fomin and Zelevinsky have shown that the coordinate rings of double Bruhat cells $G^{u,v}$ have structures of upper cluster algebras. Here, for $u,v\in W$, it is defined $G^{u,v}:=(Bu B)\cap(B_- v B_-)$. Recently, Goodearl and Yakimov have shown that the coordinate rings $\mathbb{C}[G^{u,v}]$ also have cluster algebra structures, and certain generalized minors are their initial cluster variables \cite{GY}. Generalized minors are regular functions on $G$ (see Sect.\ref{CluSect}), and we denote them by $\Delta(k;\textbf{i})$ $(k=1,2,\cdots, l(u)+l(v))$. They coincide with ordinary minors in the case $G={\rm SL}_{r+1}(\mathbb{C})$.

One purpose of our study is to reveal linkage between cluster variables of the coordinate rings $\mathbb{C}[G^{u,e}]$ and representation theories of quantum groups or quantum affine algebras. In \cite{KaN}, we gave explicit forms and path descriptions of initial cluster variables $\Delta(k;\textbf{i})$ of $\mathbb{C}[G^{u,e}]$ in the case $G={\rm SL}_{r+1}(\mathbb{C})$ (type ${\rm A}_r$ algebraic group). Using the explicit forms, we have revealed a relation between $\Delta(k;\textbf{i})$ and crystal bases in the same paper. More precisely, we found that the minors $\Delta(k;\textbf{i})$ on $G^{u,e}$ can be described as the sums of monomial realizations of certain Demazure crystals in the sense of \cite{K}. The $q$-characters of Kirillov-Reshetikhin modules of type ${\rm A}$ have similar properties to $\Delta(k;\textbf{i})$. In fact, the explicit forms of the $q$-characters were given in \cite{KNS, Nj} and we can verify they are described as the sums of monomial realizations of certain crystal bases.

The above results give us a motivation to study the case $G$ is other simple algebraic group. In \cite{KaN2}, we gave explicit formulas and path descriptions of generalized minors $\Delta(k;\textbf{i})$ on double Bruhat cells $G^{u,e}$ of type ${\rm C}_r$. 

The aim of this paper is to give explicit forms and path descriptions of $\Delta(k;\textbf{i})$ on $G^{u,e}$ of type ${\rm B}_r$. As in \cite{KaN2}, we shall only treat a Weyl group elements $u$ with the form as in \eqref{uvset} and denote its reduced word $\bfi$ by \eqref{iset}. First, we shall consider the generalized minors on reduced double Bruhat cells $L^{u,e}:=Nu N\cap B_-\subset G^{u,e}$. In \cite{B-Z}, it is shown that there exists a biregular isomorphism from
$(\mathbb{C}^{\times})^{n}$ to a Zariski open subset of $L^{u,e}$ $(n:=l(u))$ (see Theorem \ref{fp2}). We denote this isomorphism by $x^L_{\textbf{i}}$ and set $\Delta^L(k;\textbf{i}):=\Delta(k;\textbf{i})\circ x^L_{\textbf{i}}$. Second, we also define the minors $\Delta^G(k;\textbf{i})$ on double Bruhat cells $G^{u,e}$. Using a biregular isomorphism $\overline{x}^G_{\textbf{i}}$ from $H\times(\mathbb{C}^{\times})^{n}$ to a Zariski open subset of $G^{u,e}$ (Proposition \ref{gprime}), we set $\Delta^G(k;\textbf{i}):=\Delta(k;\textbf{i})\circ \overline{x}^G_{\textbf{i}}$. In Proposition \ref{gprop}, we shall show that $\Delta^G(k;\textbf{i})$ is immediately obtained from $\Delta^L(k;\textbf{i})$. Thus, we should study $\Delta^L(k;\textbf{i})$. In Proposition \ref{pathlem}, \ref{pathlem-spin}, we shall describe $\Delta^L(k;\textbf{i})$ as paths and using these path descriptions, we shall give the explicit forms of $\Delta^L(k;\textbf{i})$ in Theorem \ref{thm1}, \ref{thm2}, which are our main results. 

We will not present relations between the explicit forms of $\Delta^L(k;\textbf{i})$ and crystals here. However, we will show a relation between minors $\Delta^L(k;\textbf{i})$ and monomial realizations of crystal bases together with type ${\rm C}_r$ and ${\rm D}_r$ in forthcoming paper.

\vspace{2mm}

\nd \textbf{Acknowledgement.} I would like to express my sincere gratitude to T. Nakashima for his helpful comments and wide-ranging discussions.

\section{Fundamental representations of type ${\rm B}_r$}\label{SectFundB}

First, let us recall the fundamental representations of the complex simple Lie algebra $\ge$ of type ${\rm B}_r$ \cite{KN, N1}. We shall use them in calculations of generalized minors (see \ref{bilingen}). 
Let $I:=\{1,\cdots,r\}$ be a finite index set, $A=(a_{ij})_{i,j\in I}$ 
be the Cartan matrix of $\ge$:
\[a_{i,j}=
\begin{cases}
2 & {\rm if}\ i=j, \\
-1 & {\rm if}\ |i-j|=1\ {\rm and}\ (i,j)\neq (r,r-1), \\
-2 & {\rm if}\ (i,j)=(r,r-1), \\
0 & {\rm otherwise,}  
\end{cases}
\]
and $(\frh,\{\al_i\}_{i\in I},\{h_i\}_{i\in I})$ 
be the associated
root data 
satisfying $\al_j(h_i)=a_{ij}$ where 
$\al_i\in \frh^*$ is a simple root and 
$h_i\in \frh$ is a simple co-root. Note that $\al_i\ (i\neq r)$ are long roots and $\al_r$ is the short root. 
Let $\{\Lm_i\}_{i\in I}$ be the set of the fundamental 
weights satisfying $\Lm_i(h_j)=\del_{i,j}$, $P=\bigoplus_{i\in I}\mathbb{Z}\Lm_i$ the weight lattice and $P^*=\bigoplus_{i\in I}\mathbb{Z}h_i$ the dual weight lattice.

Define the total order on the set $J:=\{i,\ovl i|1\leq i\leq r\}\cup\{0\}$ by 
\begin{equation}\label{B-order}
 1< 2<\cd< r-1< r< 0<
 \ovl r< \ovl{r-1}< \cd< \ovl
 2< \ovl 1.
\end{equation}
For $\frg=\lan \frh,e_i,f_i(i\in I)\ran$, 
let us describe the vector representation 
$V(\Lm_1)$. Set ${\mathbf B}^{(r)}:=
\{v_i,v_{\ovl i}|i=1,2,\cd,r\}\cup\{v_0\}$ and define 
$V(\Lm_1):=\bigoplus_{v\in{\mathbf B}^{(r)}}\bbC v$. The weights of $v_i$, $v_{\ovl{i}}$ $(i=1,\cd,r)$ and $v_0$ are as follows:
\begin{equation}\label{B-wtv}
 {\rm wt}(v_i)=\Lm_i-\Lm_{i-1},\q {\rm wt}(v_{\ovl{i}})=\Lm_{i-1}-\Lm_{i}\q (1\leq i\leq r-1),
\end{equation}
\[  {\rm wt}(v_r)=2\Lm_r-\Lm_{r-1},\q {\rm wt}(v_{\ovl{r}})=\Lm_{r-1}-2\Lm_r,\q {\rm wt}(v_0)=0, \]
where $\Lm_0=0$. We define the $\frg$-action on $V(\Lm_1)$ as follows:
\begin{eqnarray}
&& h v_j=\lan h,{\rm wt}(v_j)\ran v_j\ \ (h\in P^*,\ j\in J), \\
&&f_iv_i=v_{i+1},\ f_iv_{\ovl{i+1}}=v_{\ovl i},\q
e_iv_{i+1}=v_i,\ e_iv_{\ovl i}=v_{\ovl{i+1}}
\q(1\leq i<r),\label{B-f1}\\
&&f_r v_r=v_{0},\qq 
e_r v_{\ovl r}=v_0,\qq f_r v_0=2v_{\ovl r},\qq e_r v_0=2v_r,\label{B-f2}
\end{eqnarray}
and the other actions are trivial.

Let $\Lm_i$ $(1\leq i\leq r-1)$ be the $i$-th fundamental weight of type ${\rm B}_r$.
As is well-known that the fundamental representation 
$V(\Lm_i)$ $(1\leq i\leq r-1)$
is embedded in $\wedge^i V(\Lm_1)$
with multiplicity free.
The explicit form of the highest (resp. lowest) weight 
vector $u_{\Lm_i}$ (resp. $v_{\Lm_i}$)
of $V(\Lm_i)$ is realized in 
$\wedge^i V(\Lm_1)$ as follows:
\begin{equation}
\begin{array}{ccc}\displaystyle
u_{\Lm_i}&=&v_1\wedge v_2\wedge\cdots\wedge v_i,\\
v_{\Lm_i}&=&v_{\ovl{1}}\wedge v_{\ovl{2}}\wedge\cdots \wedge v_{\ovl{i}}.
\end{array}
\label{B-h-l}
\end{equation}

The fundamental representation $V(\Lm_r)$ is called the {\it spin representation}. It can be realized as follows: Set 
\[ {\mathbf B}_{{\rm sp}}^{(r)}:=
\{(\ep_1,\cdots,\ep_r)|\ \ep_i\in\{+,- \}\q (i=1,2,\cd,r)\}, \]
\[ V_{{\rm sp}}^{(r)}:=\bigoplus_{v\in{\mathbf B}_{{\rm sp}}^{(r)}}\bbC v,\] 
and define the $\ge$-action on $V_{{\rm sp}}^{(r)}$ as follows:

\begin{equation}\label{Bsp-f0}
h_i(\ep_1,\cdots,\ep_r)=
\begin{cases}
\frac{\ep_i\cdot 1-\ep_{i+1}\cdot 1}{2}(\ep_1,\cdots,\ep_r) & {\rm if}\ i<r, \\
\ep_r(\ep_1,\cdots,\ep_r) & {\rm if}\ i=r,
\end{cases}
\end{equation}

\begin{equation}\label{Bsp-f1}
f_i(\ep_1,\cdots,\ep_r)=
\begin{cases}
(\ep_1,\cdots,\overset{i}{-},\overset{i+1}{+},\cdots,\ep_r) & {\rm if}\ \ep_i=+,\ \ep_{i+1}=-,\ i\neq r, \\
(\ep_1,\cdots,\ep_{r-1},\overset{r}{-}) & {\rm if}\ \ep_r=+,\ i=r, \\
0 & {\rm otherwise,}
\end{cases}
\end{equation}

\begin{equation}\label{Bsp-f2}
e_i(\ep_1,\cdots,\ep_r)=
\begin{cases}
(\ep_1,\cdots,\overset{i}{+},\overset{i+1}{-},\cdots,\ep_r) & {\rm if}\ \ep_i=-,\ \ep_{i+1}=+,\ i\neq r, \\
(\ep_1,\cdots,\ep_{r-1},\overset{r}{+}) & {\rm if}\ \ep_r=-,\ i=r, \\
0 & {\rm otherwise.}
\end{cases}
\end{equation}

Then the module $V_{{\rm sp}}^{(r)}$ is isomorphic to $V(\Lm_r)$ as a $\ge$-module.

\section{Factorization theorem for type ${\rm B}_r$}\label{DBCs}

In this section, we shall introduce (reduced) double Bruhat cells $G^{u,v}$, $L^{u,v}$, and their properties in the case $v=e$ and some special $u\in W$. In \cite{B-Z} and \cite{F-Z}, these properties have been proven for simply connected, connected, semisimple complex algebraic groups and arbitrary $u,v\in W$. In this paper, we treat the coordinate ring of the double Bruhat cell $G^{u,e}$, which has a structure of cluster algebra (see Sect.\ref{CluSect}). For $l\in \mathbb{Z}_{>0}$, we set $[1,l]:=\{1,2,\cdots,l\}$.

\subsection{Double Bruhat cells}\label{factpro}

Let $G$ be a simple complex algebraic group, $B$ and $B_-$ be two opposite Borel subgroups in $G$, $N\subset B$ and $N_-\subset B_-$ be their unipotent radicals, 
$H:=B\cap B_-$ a maximal torus. We set $\frg:={\rm Lie}(G)$ with the Cartan decomposition $\frg=\frn_-\oplus \frh \oplus \frn$. Let $e_i$, $f_i$ $(i\in[1,r])$ be the generators of $\frn$, $\frn_-$. For $i\in[1,r]$ and $t \in \mathbb{C}$, we set
\begin{equation}\label{xiyidef} 
x_i(t):={\rm exp}(te_i),\ \ \ y_{i}(t):={\rm exp}(tf_i).
\end{equation}

Let $W:=\lan s_i |i=1,\cdots,r \ran$ be the Weyl group of ${\rm Lie}(G)$, where
$\{s_i\}$ are the simple reflections. We identify the Weyl group $W$ with ${\rm Norm}_G(H)/H$. An element 
\begin{equation}\label{smpl}
\ovl{s_i}:=x_i(-1)y_i(1)x_i(-1)
\end{equation}
is in ${\rm Norm}_G(H)$, which is representative of $s_i\in W={\rm Norm}_G(H)/H$ \cite{N1}. For a reduced expression $w=s_{i_1}\cdots s_{i_n}\in W$, we define $\ovl{w}:=\ovl{s_{i_1}}\cdots\ovl{s_{i_n}}$. We call $l(w):=n$ the length of $w$.

We have two kinds of Bruhat decompositions of $G$ as follows:
\[ G=\displaystyle\coprod_{u \in W}B\ovl{u}B=\displaystyle\coprod_{u \in W}B_-\ovl{u}B_- .\]
Then, for $u$, $v\in W$, 
we define the {\it double Bruhat cell} $G^{u,v}$ as follows:
\[ G^{u,v}:=B\ovl{u}B \cap B_-\ovl{v}B_-. \]
This is biregularly isomorphic to a Zariski open subset of 
an affine space of dimension $r+l(u)+l(v)$ \cite[Theorem 1.1]{F-Z}.

We also define the {\it reduced double Bruhat cell} $L^{u,v}$ as follows:
\[ L^{u,v}:=N\ovl{u}N \cap B_-\ovl{v}B_- \subset G^{u,v}. \] 
As is the case with $G^{u,v}$, $L^{u,v}$ is 
biregularly isomorphic to a Zariski open subset of an 
affine space of dimension $l(u)+l(v)$ \cite[Proposition 4.4]{B-Z}.

\begin{defn}\label{redworddef}
Let $u=s_{i_1}\cdots s_{i_n}$ be a reduced expression of $u\in W$ $(i_1,\cdots,i_n\in [1,r])$. Then the finite sequence 
\[ \textbf{i}:=(i_1,\cdots,i_n) \]
is called a {\it reduced word} for $u$.
\end{defn}

For example, the sequence $(1,2,3,1,2,3,1,2,3)$ is a reduced word of the longest element $s_1s_2s_3s_1s_2s_3s_1s_2s_3$ of the Weyl group of type ${\rm B}_3$. For the case of type ${\rm B}_r$, we fix the reduced word $\textbf{i}_0$ of the longest element as follows:

\begin{equation}\label{redwords}
\textbf{i}_0=(1,2,\cdots,r-1,r)^r 
\end{equation}
In this paper, we mainly treat (reduced) Double Bruhat cells of the form $G^{u,e}:=B\ovl{u}B \cap B_-$, $L^{u,e}:=N\ovl{u}N \cap B_-$, and the element $u\in W$ whose reduced word can be written as a left factor of $\textbf{i}_0$.

\subsection{Factorization theorem for type ${\rm B}_r$}\label{factpro0}

In this subsection, we consider the case of type ${\rm B}_r$ ($G={\rm SO}_{2r+1}(\mathbb{C})$) and introduce the isomorphisms between the double Bruhat cell $G^{u,e}$ and $H\times (\mathbb{C}^{\times})^{l(u)}$, and between $L^{u,e}$ and $(\mathbb{C}^{\times})^{l(u)}$. 

For a reduced word $\textbf{i}=(i_1, \cdots ,i_n)$ of $u$
($i_1,\cdots,i_n\in[1,r]$), 
we define a map $x^G_{\textbf{i}}:H\times \mathbb{C}^n \rightarrow G$ as 
\begin{equation}\label{xgdef}
x^G_{\textbf{i}}(a; t_1, \cdots, t_n):=a\cdot y_{i_1}(t_1)\cdots y_{i_n}(t_n).
\end{equation}

\begin{thm}\label{fp}${\cite[Theorem\ 1.2]{F-Z}}$ Let $u\in W$ be an element whose reduced word ${\rm \bf{i}}$ can be written as a left factor of ${\rm \bf{i}}_0$ $(\ref{redwords})$. The map $x^G_{{\rm \bf{i}}}$ defined above can be restricted to a biregular isomorphism between $H\times (\mathbb{C}^{\times})^{l(u)}$ and a Zariski open subset of $G^{u,e}$. 
\end{thm}

Next, for $i \in [1,r]$ and $t\in \mathbb{C}^{\times}$, we define as follows:
\begin{equation}\label{alxmdef}
\alpha_i^{\vee}(t):=t^{h_i},\ \ 
 x_{-i}(t):=y_{i}(t)\alpha_i^{\vee}(t^{-1}).
\end{equation}
For $\textbf{i}=(i_1, \cdots ,i_n)$
($i_1,\cdots,i_n\in[1,r]$), 
we define a map $x^L_{\textbf{i}}:\mathbb{C}^n \rightarrow G$ as 
\begin{equation}\label{xldef}
x^L_{\textbf{i}}(t_1, \cdots, t_n):=x_{-i_1}(t_1)\cdots x_{-i_n}(t_n).
\end{equation}
We have the following theorem which is similar to the previous one.
\begin{thm}\label{fp2}${\cite[Proposition\ 4.5]{B-Z}}$
Let $u\in W$ be an element whose reduced word ${\rm \bf{i}}$ can be written as a left factor of ${\rm \bf{i}}_0$ $(\ref{redwords})$.
The map $x^L_{{\rm \bf{i}}}$ defined above can be restricted to a
biregular isomorphism between $ (\mathbb{C}^{\times})^{l(u)}$ 
and a Zariski open subset of $L^{u,e}$. 
\end{thm}

We define a map
$\ovl{x}^G_{\textbf{i}}:H\times(\mathbb{C}^{\times})^{n}\rightarrow
G^{u,e}$ as
\[ \ovl{x}^G_{\textbf{i}}(a;t_1,\cdots,t_n)
=ax^L_{\textbf{i}}(t_1,\cdots,t_n), \]
where $a\in H$ and $(t_1,\cdots,t_n)\in (\mathbb{C}^{\times})^{n}$.
\begin{prop}\label{gprime}
In the above setting, the map $\ovl{x}^G_{{\rm \bf{i}}}$ is a biregular isomorphism between $H\times(\mathbb{C}^{\times})^{n}$ and a Zariski open subset of $G^{u,e}$.
\end{prop}

\nd
{\sl Proof.}

In this proof, we use the notation
\[ (Y_{1,1},\cdots,Y_{1,r},\cdots,Y_{m-1,1},\cdots,Y_{m-1,r},Y_{m,1},\cdots,Y_{m,i_n})\in (\mathbb{C}^{\times})^{n} \]
for variables instead of $(t_1,\cdots,t_n)$.

We define a map
$\phi:H\times(\mathbb{C}^{\times})^{n}\rightarrow
H\times(\mathbb{C}^{\times})^{n}$ as follows: For 
\[
\textbf{Y}:=(a;Y_{1,1},\cdots,Y_{1,r}, \cdots,Y_{m,1},\cdots,Y_{m,i_n}),
\]
we define
$\phi(\textbf{Y})=(\Phi_a(\textbf{Y});\Phi_{1,1}(\textbf{Y}),\cdots,\Phi_{1,r}(\textbf{Y}),\cdots,\Phi_{m,1}(\textbf{Y}),\cdots,\Phi_{m,i_n}(\textbf{Y}))$ as
\[ \Phi_a(\textbf{Y}):=a\cdot\left(\prod^{m-1}_{j=1}{\al_1^{\vee}(Y_{j,1})^{-1}\cdots
\al_r^{\vee}(Y_{j,r})^{-1}}\right)\cdot\al_1^{\vee}(Y_{m,1})^{-1}\cdots
\al_{i_n}^{\vee}(Y_{m,i_n})^{-1}, \]
and for $1\leq s\leq m$,
\begin{equation}\label{mbase0} 
\Phi_{s,l}(\textbf{Y}):=
\begin{cases}
\frac{(Y_{s+1,l-1}Y_{s+2,l-1}\cdots Y_{m,l-1})(Y_{s,l+1}Y_{s+1,l+1}\cdots Y_{m,l+1})}{Y_{s,l}(Y_{s+1,l}\cdots Y_{m,l})^{2}} & {\rm if}\ 1\leq l\leq r-2, \\
\frac{(Y_{s+1,r-2}Y_{s+2,r-2}\cdots Y_{m,r-2})(Y_{s,r}Y_{s+1,r}\cdots Y_{m,r})^{2}}{Y_{s,r-1}(Y_{s+1,r-1}\cdots Y_{m,r-1})^{2}} & {\rm if}\ l=r-1,\\
\frac{(Y_{s+1,r-1}Y_{s+2,r-1}\cdots Y_{m,r-1})}{Y_{s,r}(Y_{s+1,r}\cdots Y_{m,r})^{2}} & {\rm if}\ l=r, \\
\end{cases}
\end{equation}
where in (\ref{mbase0}), if we see the variables $Y_{\zeta,0}$ $(1\leq\zeta\leq m)$ and $Y_{m,\xi}$ $(i_n<\xi)$, then we understand $Y_{\zeta,0}=Y_{m,\xi}=1$. For example, $Y_{s+1,l-1}=1$ in the case $l=1$.
Note that $\phi$ is a biregular isomorphism since we can recurrently construct the inverse map $\psi:H\times(\mathbb{C}^{\times})^{n}\rightarrow
H\times(\mathbb{C}^{\times})^{n}$, $\textbf{Y}\mapsto (\Psi_a(\textbf{Y});\Psi_{1,1}(\textbf{Y}),\cdots,\Psi_{m,i_n}(\textbf{Y}))$ of $\phi$ as follows: The definition (\ref{mbase0}) implies that $\Phi_{m,i_n}(\textbf{Y})=\frac{1}{Y_{m,i_n}}$, and hence $Y_{m,i_n}=\frac{1}{\Psi_{m,i_n}(\textbf{Y})}$. So we set $\Psi_{m,i_n}(\textbf{Y})=\frac{1}{Y_{m,i_n}}$. Suppose that we can construct $\Psi_{m,i_n}(\textbf{Y}),\Psi_{m,i_n-1}(\textbf{Y}),\cdots$,
$\Psi_{m,1}(\textbf{Y}),\cdots\Psi_{s+1,r}(\textbf{Y})$,
$\cdots,\Psi_{s+1,1}(\textbf{Y}),\Psi_{s,r}(\textbf{Y}),\cdots$,
$\Psi_{s,l+1}(\textbf{Y})$. Then we define
\[ \Psi_{s,l}(\textbf{Y}):=
\begin{cases}
\frac{(\Psi_{s+1,l}(\textbf{Y})\cdots \Psi_{m,l}(\textbf{Y}))^{2}}
{Y_{s,l}(\Psi_{s+1,l-1}(\textbf{Y})\Psi_{s+2,l-1}(\textbf{Y})\cdots \Psi_{m,l-1}(\textbf{Y}))(\Psi_{s,l+1}(\textbf{Y})\cdots \Psi_{m,l+1}(\textbf{Y}))}
 & {\rm if}\ 1\leq l\leq r-2, \\
\frac{(\Psi_{s+1,r-1}(\textbf{Y})\cdots \Psi_{m,r-1}(\textbf{Y}))^{2}}
{Y_{s,r-1}(\Psi_{s+1,r-2}(\textbf{Y})\Psi_{s+2,r-2}(\textbf{Y})\cdots \Psi_{m,r-2}(\textbf{Y}))(\Psi_{s,r}(\textbf{Y})\cdots \Psi_{m,r}(\textbf{Y}))^2}
 & {\rm if}\ l=r-1, \\
\frac{(\Psi_{s+1,r}(\textbf{Y})\cdots \Psi_{m,r}(\textbf{Y}))^{2}}
{Y_{s,r}(\Psi_{s+1,r-1}(\textbf{Y})\Psi_{s+2,r-1}(\textbf{Y})\cdots \Psi_{m,r-1}(\textbf{Y}))}
 & {\rm if}\ l=r.
\end{cases}
\]
We also define
\[ \Psi_{a}(\textbf{Y}):=
a\cdot\left(\prod^{m-1}_{j=1}{\al_1^{\vee}(\Psi_{j,1}(\textbf{Y}))\cdots
\al_r^{\vee}(\Psi_{j,r}(\textbf{Y}))}\right)\cdot\al_1^{\vee}(\Psi_{m,1}(\textbf{Y}))\cdots
\al_{i_n}^{\vee}(\Psi_{m,i_n}(\textbf{Y})).
\]
Then, we get the inverse map $\psi$ of $\phi$.

Let us prove
\[ \ovl{x}^G_{\textbf{i}}(\textbf{Y})=(x^G_{\textbf{i}}\circ\phi)(\textbf{Y}), \]
which implies that $\ovl{x}^G_{\textbf{i}}:H\times(\mathbb{C}^{\times})^{n}\rightarrow G^{u,e}$ is biregular isomorphism by Theorem \ref{fp}.

First, it is known that
\begin{equation}\label{base2}
\al_i^{\vee}(c)^{-1}y_{j}(t)=\begin{cases}
	y_{i}(c^2t)\al_i^{\vee}(c)^{-1} & {\rm if}\ i=j, \\
	y_{j}(c^{-1}t)\al_i^{\vee}(c)^{-1} & {\rm if}\ |i-j|=1\ {\rm and}\ (i,j)\neq (r,r-1), \\
	y_{j}(c^{-2}t)\al_i^{\vee}(c)^{-1} & {\rm if}\ (i,j)= (r,r-1), \\
	y_{j}(t)\al_i^{\vee}(c)^{-1} & {\rm otherwise},
\end{cases}
\end{equation}
for $1\leq i,\ j\leq r$ and $c,\ t\in \mathbb{C}^{\times}$.

On the other hand, it follows from the definition (\ref{xgdef}) of $x^G_{\textbf{i}}$ and $(\ref{mbase0})$ that
\begin{multline}\label{xigp}(x^G_{\textbf{i}}\circ\phi)(\textbf{Y})\\
=a\times\left(\prod^{m-1}_{j=1}{\al_1^{\vee}(Y_{j,1})^{-1}\cdots
\al_r^{\vee}(Y_{j,r})^{-1}}\right)\cdot\al_1^{\vee}(Y_{m,1})^{-1}\cdots
\al_{i_n}^{\vee}(Y_{m,i_n})^{-1}
\\ 
 \times y_{1}(\Phi_{1,1}(\textbf{Y}))y_{2}(\Phi_{1,2}(\textbf{Y}))\cdots
 y_{r}(\Phi_{1,r}(\textbf{Y}))\cdots y_{1}(\Phi_{m,1}(\textbf{Y}))
\cdots y_{i_n}(\Phi_{m,i_n}(\textbf{Y})).
\end{multline}

For each $s$ and $l$ $(1\leq s\leq m,\ 1\leq l\leq r)$, we can move 
\begin{multline*}
\al_{l}^{\vee}(Y_{s,l})^{-1}
\al_{l+1}^{\vee}(Y_{s,l+1})^{-1}
\cdots\al_{r}^{\vee}(Y_{s,r})^{-1}\\
\cdot\left(\prod^{m-1}_{j=s+1}\al_1^{\vee}(Y_{j,1})^{-1}\cdots \al_r^{\vee}(Y_{j,r})^{-1}\right)\cdot
\al_1^{\vee}(Y_{m,1})^{-1}\cdots\al_{i_n}^{\vee}(Y_{m,i_n})^{-1}
\end{multline*}
to the right of
$y_{l}(\Phi_{s,l}(\textbf{Y}))$ by using the relations (\ref{base2}). 
For example,
\begin{flushleft}
$\al_1^{\vee}(Y_{m,1})^{-1}\cdots
\al_{i_n}^{\vee}(Y_{m,i_n})^{-1} y_{l}(\Phi_{s,l}(\textbf{Y}))=$
\end{flushleft}
\[
\begin{cases}
y_{l}\left(\frac{Y_{m,l}^2}{Y_{m,l-1}Y_{m,l+1}}\Phi_{s,l}(\textbf{Y})\right)\al_1^{\vee}(Y_{m,1})^{-1}\cdots
\al_{i_n}^{\vee}(Y_{m,i_n})^{-1} &\ {\rm if}\ 1\leq l\leq r-2,\\
y_{r-1}\left(\frac{Y_{m,r-1}^2}{Y_{m,r-2}Y_{m,r}^2}\Phi_{s,r}(\textbf{Y})\right)\al_1^{\vee}(Y_{m,1})^{-1}\cdots
\al_{i_n}^{\vee}(Y_{m,i_n})^{-1} &\ {\rm if}\ l=r-1,\\
y_{l}\left(\frac{Y_{m,r}^2}{Y_{m,r-1}}\Phi_{s,l}(\textbf{Y})\right)\al_1^{\vee}(Y_{m,1})^{-1}\cdots
\al_{i_n}^{\vee}(Y_{m,i_n})^{-1} &\ {\rm if}\ l=r.
\end{cases}
\]

Repeating this argument, in the case $1\leq l\leq r-2$ or $l=r$, we have
\begin{multline*}
\al_{l}^{\vee}(Y_{s,l})^{-1}
\al_{l+1}^{\vee}(Y_{s,l+1})^{-1}
\cdots\al_{r}^{\vee}(Y_{s,r})^{-1}\\
\times\left(\prod^{m-1}_{j=s+1}\al_1^{\vee}(Y_{j,1})^{-1}\cdots \al_r^{\vee}(Y_{j,r})^{-1}\right)\cdot
\al_1^{\vee}(Y_{m,1})^{-1}\cdots\al_{i_n}^{\vee}(Y_{m,i_n})^{-1}y_{l}(\Phi_{s,l}(\textbf{Y}))\\
=y_{l}\left(\frac{(Y_{s,l}Y_{s+1,l}\cdots Y_{m-1,l}Y_{m,l})^2}{(Y_{s+1,l-1}\cdots Y_{m-1,l-1}Y_{m,l-1})(Y_{s,l+1}\cdots Y_{m-1,l+1}Y_{m,l+1})}\Phi_{s,l}(\textbf{Y})\right)\cdot \al_{l}^{\vee}(Y_{s,l})^{-1}\\
\times\al_{l+1}^{\vee}(Y_{s,l+1})^{-1}
\cdots\al_{r}^{\vee}(Y_{s,r})^{-1}
\cdot\left(\prod^{m-1}_{j=s+1}\al_1^{\vee}(Y_{j,1})^{-1}\cdots \al_r^{\vee}(Y_{j,r})^{-1}\right)\cdot
\al_1^{\vee}(Y_{m,1})^{-1}\cdots\al_{i_n}^{\vee}(Y_{m,i_n})^{-1}.
\end{multline*}

Note that $\frac{(Y_{s,l}Y_{s+1,l}\cdots Y_{m-1,l}Y_{m,l})^2}{(Y_{s+1,l-1}\cdots Y_{m-1,l-1}Y_{m,l-1})(Y_{s,l+1}\cdots Y_{m-1,l+1}Y_{m,l+1})}\Phi_{s,l}(\textbf{Y})=Y_{s,l}$, which implies 
\begin{equation}\label{tradeeq}
\al_{l}^{\vee}(Y_{s,l})^{-1}\cdots \al_{i_n}^{\vee}(Y_{m,i_n})^{-1}y_{l}(\Phi_{s,l}(\textbf{Y}))=y_{l}(Y_{s,l})\al_{l}^{\vee}(Y_{s,l})^{-1}\cdots \al_{i_n}^{\vee}(Y_{m,i_n})^{-1}.
\end{equation}

In the case $l=r-1$, we can also verify the relation (\ref{tradeeq}) similarly. Thus, by (\ref{xigp}) and (\ref{tradeeq}), we have
\begin{eqnarray*}
&&(x^G_{\textbf{i}}\circ\phi)(\textbf{Y})=a\cdot
y_{1}(Y_{1,1})\al_1^{\vee}(Y_{1,1})^{-1}\cdots y_{r}(Y_{1,r})\al_r^{\vee}(Y_{1,r})^{-1}\times\cdots \\
&&\qq\times y_{1}(Y_{m,1})\al_1^{\vee}(Y_{m,1})^{-1}\cdots
 y_{i_n}(Y_{m,i_n})\al_{i_n}^{\vee}(Y_{m,i_n})^{-1}\\
&&
=a\cdot x_{-1}(Y_{1,1})\cdots x_{-r}(Y_{1,r})\cdots x_{-1}(Y_{m,1})\cdots
x_{-i_n}(Y_{m,i_n})\\
&&=\ovl{x}^G_{\textbf{i}}(\textbf{Y}).\qq\qq\qq\qq\qq\qq\qq\qq\qq \qed
\end{eqnarray*}

\section{Cluster algebras and generalized minors}\label{CluSect}
Following \cite{A-F-Z,F-Z,FZ2,M-M-A}, we review the definitions of (upper) cluster algebras and their generators called cluster variables. It is known that the coordinate rings of double Bruhat cells have upper cluster algebra structures, and generalized minors are their cluster variables \cite{A-F-Z}. We will refer to explicit forms of certain cluster variables on double Bruhat cells in Sect.\ref{gmc}.

We set $[1,l]:=\{1,2,\cdots,l\}$ and $[-1,-l]:=\{-1,-2,\cdots,-l\}$ for $l\in \mathbb{Z}_{>0}$. For $n,m\in \mathbb{Z}_{>0}$, let $x_1, \cdots ,x_n,x_{n+1}, \cdots
,x_{n+m}$ be commuting variables and $\mathcal{P}$ be a free multiplicative
abelian group generated by $x_{n+1},\cdots,x_{n+m}$. We set ${\mathbb
Z}\cP:={\mathbb Z}[x_{n+1}^{\pm1}, \cdots ,x_{n+m}^{\pm1}]$. Let $\cF:=\mathbb{C}(x_{1}, \cdots ,x_{n},x_{n+1},\cdots,x_{n+m})$ 
be the field of rational functions.

\subsection{Cluster algebras of geometric type}

In this subsection, we recall the definitions of (upper) cluster algebras. Let $\tilde{B}=(b_{ij})_{1\leq i\leq n+m,\ 1\leq j \leq n}$ be an $(n+m)\times
n$ integer matrix. The {\it principal part} $B$ of $\tilde{B}$ is obtained from $\tilde{B}$ by deleting the last $m$ rows. For $\tilde{B}$ and $k\in [1,n]$, the new $(n+m)\times n$ integer matrix $\mu_k(\tilde{B})=(b'_{ij})$ is defined by
\[b_{ij}':=
\begin{cases}
	-b_{ij} & {\rm if}\ i=k\ {\rm or}\ j=k, \\
	b_{ij}+\frac{|b_{ik}|b_{kj}+b_{ik}|b_{kj}|}{2} & {\rm otherwise}.
\end{cases}
\]
One calls $\mu_k(\tilde{B})$ the {\it matrix mutation} in direction $k$ of $\tilde{B}$. If there exists a positive 
integer diagonal matrix $D$ such that $DB$ is skew symmetric, we say $B$ is {\it skew symmetrizable}. It is easily verified that if $\tilde{B}$ has a skew symmetrizable principal part then $\mu_k(\tilde{B})$ also has a skew symmetrizable principal part${\cite[Proposition\ 3.6]{M-M-A}}$. We can also verify that $\mu_k\mu_k(\tilde{B})=\tilde{B}$. Define $\textbf{x}:=(x_1,\cdots,x_{n+m})$ and we call the pair $(\textbf{x}, \tilde{B})$ {\it initial seed}. For $1\leq k\leq n$, a new cluster variable $x_k'$ is defined by
\[ x_k x_k' = 
\prod_{1\leq i \leq n+m,\ b_{ik}>0} x_i^{b_{ik}}
+\prod_{1\leq i \leq n+m,\ b_{ik}<0} x_i^{-b_{ik}}. \]
Let $\mu_k(\textbf{x})$ be the set of variables obtained from $\textbf{x}$ by replacing $x_k$ by $x'_k$. Ones call the pair $(\mu_k(\textbf{x}), \mu_k(\tilde{B}))$ the {\it mutation} in direction $k$ of the seed $(\textbf{x}, \tilde{B})$.

Now, we can repeat this process of mutation and obtain a set of seeds inductively. Hence, each seed consists of an $(n+m)$-tuple of variables and a matrix. Ones call this $(n+m)$-tuple and matrix {\it cluster} and {\it exchange matrix} respectively. Variables in the clusters are called {\it cluster variables}.

\begin{defn}${\cite{F-Z, M-M-A}}$\label{clusterdef}
Let $\tilde{B}$ be a integer matrix whose principal part is skew symmetrizable and $\Sigma=(\textbf{x},\tilde{B})$ a seed. We set ${\mathbb A}:={\mathbb Z}\cP$. The cluster algebra (of geometric type) $\cA=\cA(\Sigma)$ over $\mathbb A$ associated with
seed $\Sigma$ is defined as the ${\mathbb A}$-subalgebra of $\cF$ generated by all cluster variables in all seeds which can be obtained from $\Sigma$ by sequences of mutations.
\end{defn}

For a seed $\Sigma$, we define ${\mathbb Z}\cP$-subalgebra $\UU(\Sigma)$ of $\cF$ by
\[ \UU(\Sigma):={\mathbb Z}\cP[\textbf{x}^{\pm 1}] \cap {\mathbb Z}\cP[\textbf{x}_1^{\pm 1}] \cap \cdots \cap {\mathbb Z}\cP[\textbf{x}_n^{\pm 1}] .\]
 Here, ${\mathbb Z}\cP[\textbf{x}^{\pm 1}]$ is the Laurent polynomial ring in \textbf{x} and $\textbf{x}_k:=\mu_k(\textbf{x})$.

\begin{defn}${\cite{F-Z, M-M-A}}$\label{upper}
Let $\Sigma$ be the seed in Definition \ref{clusterdef}. We define an {\it upper cluster algebra} $\ovl{\cA}=\ovl{\cA}(\Sigma)$ as the intersection of the subalgebras $\UU(\Sigma')$ for all seeds $\Sigma'$ which can be obtained from $\Sigma$ by sequences of mutations.
\end{defn}
Following the inclusion relation holds \cite{A-F-Z}:
\[ \cA(\Sigma) \subset \ovl{\cA}(\Sigma) .\]

\subsection{Upper cluster algebra $\ovl{\cA}(\textbf{i})$}

Let $G$ be a simple classical algebraic group, $\ge:={\rm Lie}(G)$ and $A=(a_{i,j})$ be its Cartan matrix. In Definition \ref{redworddef}, we define a reduced word $\textbf{i}=(i_1,\cdots,i_{l(u)})$ for an element $u$ of Weyl group $W$. In this subsection, we define the upper cluster algebra $\ovl{\cA}(\textbf{i})$, which obtained from $\textbf{i}$. It satisfies that $\ovl{\cA}(\textbf{i})\otimes \mathbb{C}$ is isomorphic to the coordinate ring $\mathbb{C}[G^{u,e}]$ of the double Bruhat cell \cite{A-F-Z}. Let $i_k$ $(k\in[1,l(u)])$ be the $k$-th index of $\textbf{i}$ from the left. For $t\in [-1,-r]$, we set $i_t:=t$.

For $k\in[-1,-r]\cup[1,l(u)]$, we denote by $k^+$ the smallest index $l$ such that $k<l$ and $|i_l|=|i_k|$. For example, if $\textbf{i}=(1,2,3,1,2)$ then, $1^+=4$, $2^+=5$ and $3^+$ is not defined. We define a set e(\textbf{i}) as 
\[ e(\textbf{i}):= \{k\in[1,l(u)]| k^+\ {\rm is\ well}-{\rm defined}\}.  \]
Following \cite{A-F-Z}, we define a quiver $\Gamma_{\textbf{i}}$ as follows. The vertices of $\Gamma_{\textbf{i}}$ are the numbers $[-1,-r]\cup[1,l(u)]$. For two vertices $k\in [-1,-r]\cup[1,l(u)]$ and $l\in[1,l(u)]$ with $k<l$, there exists an arrow $k\rightarrow l$ (resp. $l\rightarrow k$) if and only if $l=k^+$ (resp. $l<k^+<l^+$ and $a_{i_k,i_l}<0$). Next, let us define a matrix $\tilde{B}=\tilde{B}(\textbf{i})$. 
\begin{defn}
Let $\tilde{B}(\textbf{i})$ be an integer matrix with rows labelled by all the indices in $[-1,-r]\cup [1,l(u)]$ and columns labelled by all the indices in $e(\textbf{i})$. For $k\in[-1,-r]\cup [1,l(u)]$ and $l\in e(\textbf{i})$, an entry $b_{kl}$ of $\tilde{B}(\textbf{i})$ is determined as follows: If there exists an arrow $k\rightarrow l$ (resp. $l\rightarrow k$) in $\Gamma_{\textbf{i}}$, then
\[
b_{kl}:=\begin{cases}
		1\ ({\rm resp.}\ -1)& {\rm if}\ |i_k|=|i_l|, \\
		-a_{|i_k||i_l|}\ ({\rm resp.}\ a_{|i_k||i_l|})& {\rm if}\ |i_k|\neq|i_l|.
	\end{cases}
\]
If there exist no arrows between $k$ and $l$, we set $b_{kl}=0$.
\end{defn}

\begin{prop}\label{propss}${\cite[Proposition\ 2.6]{A-F-Z}}$
$\tilde{B}({\rm \bf{i}})$ is skew symmetrizable. 
\end{prop}

By Definition $\ref{upper}$ and Proposition \ref{propss}, we can construct the upper cluster algebra:

\begin{defn}
We denote this upper cluster algebra by $\ovl{\cA}(\textbf{i})$.
\end{defn}

\subsection{Generalized minors and bilinear form}\label{bilingen}

Set $\bar{\cA}(\textbf{i})_{\mathbb{C}}:=\bar{\cA}(\textbf{i})\otimes \mathbb{C}$. It is known that the coordinate ring $\mathbb{C}[G^{u,e}]$ of the double Bruhat cell is isomorphic to $\bar{\cA}(\textbf{i})_{\mathbb{C}}$ (Theorem \ref{clmainthm}). To describe this isomorphism explicitly, we need generalized minors.  

We set $G_0:=N_-HN$, and let $x=[x]_-[x]_0[x]_+$ with $[x]_-\in N_-$, $[x]_0\in H$, $[x]_+\in N$ be the corresponding decomposition. 

\begin{defn}
For $i\in[1,r]$ and $u'\in W$, the generalized minor $\Delta_{u\Lambda_i,\Lambda_i}$ is the regular function on $G$ whose restriction to the open set $\ovl{u} G_0$ is given by $\Delta_{u\Lambda_i,\Lambda_i}(x)=([\ovl{u}^{-1}x ]_0)^{\Lambda_i}$. Here, $\Lambda_i$ is the $i$-th  fundamental weight. In particular, we write $\Delta_{\Lambda_i}:=\Delta_{\Lambda_i,\Lambda_i}$ and call it a {\it principal minor}.
\end{defn}

We can calculate generalized minors by using a bilinear form in the fundamental representation of $\ge={\rm Lie}(G)$ (see Sect.\ref{SectFundB}). Let $\omega:\ge\to\ge$ be the anti involution 
\[
\omega(e_i)=f_i,\q
\omega(f_i)=e_i,\q \omega(h)=h,
\] and extend it to $G$ by setting
$\omega(x_i(c))=y_{i}(c)$, $\omega(y_{i}(c))=x_i(c)$ and $\omega(t)=t$
$(t\in H)$. Here, $x_i$ and $y_i$ were defined in \ref{factpro} (\ref{xiyidef}). There exists a $\ge$ (or $G$)-invariant bilinear form on the
fundamental representation $V(\Lm_i)$ of $\ge$ such that 
\[
 \lan au,v\ran=\lan u,\omega(a)v\ran,
\q\q(u,v\in V(\lm),\,\, a\in \ge\ (\text{or }G)).
\]
For $g\in G$, 
we have the following simple fact:
\[
 \Del_{\Lm_i}(g)=\lan gu_{\Lm_i},u_{\Lm_i}\ran,
\]
where $u_{\Lm_i}$ is a properly normalized highest weight vector in
$V(\Lm_i)$. Hence, for $w\in W$, we have
\begin{equation}\label{minor-bilin}
 \Del_{w\Lm_i,\Lm_i}(g)=
\Del_{\Lm_i}({\ovl w}^{-1}g)=
\lan {\ovl w}^{-1}g\cdot u_{\Lm_i},u_{\Lm_i}\ran
=\lan g\cdot u_{\Lm_i}\, ,\, \ovl{w}\cdot u_{\Lm_i}\ran,
\end{equation}
where $\ovl w$ is the one we defined in \ref{factpro} (\ref{smpl}), and note that $\omega(\ovl s_i^{\pm})=\ovl s_i^{\mp}$.

\subsection{Cluster algebras on Double Bruhat cells}\label{cdb}

For $u=s_{i_1}s_{i_2}\cdots s_{i_n}$ and $k\in [1,l(u)]$, we set
\begin{equation}\label{inc}
u_{\leq k}=u_{\leq k}(\textbf{i}):=s_{i_1}s_{i_2}\cdots s_{i_k}.
\end{equation}
For $k \in [-1,-r]$, we set $u_{\leq k}:=e$ and $i_k:=k$. For $k \in [-1,-r]\cup [1,l(u)]$, we define
\[ \Delta(k;\textbf{i})(x):=\Delta_{u_{\leq k} \Lambda_{|i_k|},\Lambda_{|i_k|}}(x).  \]

We set
\[ F(\textbf{i}):=\{ \Delta(k;\textbf{i})(x)|k \in [-1,-r]\cup[1,l(u)] \}. \]
It is known that the set $F(\textbf{i})$ is an algebraically independent generating set for the field of rational functions $\mathbb{C}(G^{u,e})$ \cite[Theorem 1.12]{F-Z}. 

\begin{thm}\label{clmainthm}${\cite[Theorem\ 2.10]{A-F-Z}}$
The isomorphism of fields $\varphi :F \rightarrow \mathbb{C}(G^{u,e})$ defined by $\varphi (x_k)=\Delta(k;{\rm \bf{i}})\ (k \in [-1,-r]\cup [1,l(u)] )$ restricts to an isomorphism of algebras $\bar{\cA}({\rm \bf{i}})_{\mathbb{C}}\rightarrow \mathbb{C}[G^{u,e}]$.
\end{thm}

\begin{rem}\label{upclu}
In $\cite{GY}$, it is shown that $\mathbb{C}[G^{u,v}]$ has a structure of cluster algebra not only upper cluster algebra. 
\end{rem}

\section{Explicit forms of cluster variables on D.B cells of type ${\rm B}_r$}\label{gmc}

In the rest of the paper, we consider the case
$G={\rm SO}_{2r+1}(\mathbb{C})$. 
Let $u\in W$ be
\begin{equation}\label{uvset}
u=(s_1s_2\cdots s_r)^{m-1}s_1\cdots s_{i_n},
\end{equation} 
where $n=l(u)$, $1\leq i_n\leq r$ and $1\leq m\leq r$. Let
\begin{equation}\label{iset}
 \textbf{i}=(\underbrace{1,\cdots,r}_{1{\rm \,st\
  cycle}},\underbrace{1,\cdots,r}_{2{\rm \,nd\
  cycle}},\cdots,\underbrace{1,\cdots,r}_{m-1{\rm \,th\ cycle}},
\underbrace{1,\cdots,i_n}_{m{\rm \,th\ cycle}})
\end{equation}
be a reduced word for $u$, that is, \textbf{i} is the left factor of $(1,2,3,\cdots,r)^r$.  Let $i_k$ be the $k$-th index of  $\textbf{i}$ from the left, and belong to $m'$-th cycle $(m'\leq m)$. As we shall show in Lemma \ref{gmlem}, we may assume $i_n=i_k$. 

By Theorem \ref{clmainthm} and Remark \ref{upclu}, the coordinate ring $\mathbb{C}[G^{u,e}]$ has the structure of (upper) cluster algebra and $\{\Delta(k;\textbf{i})\}$ are its initial cluster
variables. Each $\Delta(k;\textbf{i})$ is a regular function on
$G^{u,e}$. On the other hand, by Proposition \ref{gprime} (resp. Theorem
\ref{fp2}), we can consider $\Delta(k;\textbf{i})$ as a function on
$H\times
(\mathbb{C}^{\times})^{l(u)}$ (resp. $(\mathbb{C}^{\times})^{l(u)}$). Here,
we change the variables of 
$\{\Delta(k;\textbf{i})\}$ as follows: 
\begin{defn}\label{gendef}For $a\in H$ and
\begin{multline}\label{yset}
 \textbf{Y}:=(Y_{1,1},Y_{1,2},\cdots,Y_{1,r},Y_{2,1},Y_{2,2},\cdots,Y_{2,r}, \\
\cdots,Y_{m-1,1},\cdots,Y_{m-1,r},Y_{m,1},\cdots,Y_{m,i_n})\in  (\mathbb{C}^{\times})^{n}, 
\end{multline}
we set 
\[ 
\Delta^G(k;\textbf{i})(a,\textbf{Y}):=(\Delta(k;\textbf{i})\circ
 \ovl{x}^G_{\textbf{i}})(a,\textbf{Y}),\ \ 
\Delta^L(k;\textbf{i})(\textbf{Y}):=(\Delta(k;\textbf{i})\circ
 x^L_{\textbf{i}})(\textbf{Y}),
\]
where $\ovl{x}^G_{\textbf{i}}$ and $x^L_{\textbf{i}}$ are as in \ref{factpro0}.
\end{defn}
We will describe the function $\Delta^L(k;\textbf{i})(\textbf{Y})$ explicitly since $\Delta^G(k;\textbf{i})(a,\textbf{Y})$ is immediately obtained from $\Delta^L(k;\textbf{i})(\textbf{Y})$ (Proposition \ref{gprop}).

\begin{rem}\label{importantrem}
If we see the variables $Y_{s,0}$, $Y_{s,r+1}$ 
$(1\leq s\leq m)$ then 
we understand $Y_{s,0}=Y_{s,r+1}=1$. For example, if $i=1$ then $Y_{s,i-1}=1$.
\end{rem}

\subsection{Generalized minors $\Delta^G(k;\textbf{i})(a,\textbf{Y})$ and $\Delta^L(k;\textbf{i})(\textbf{Y})$}

In this subsection, we shall prove some properties of the minors $\Delta^G$ and $\Delta^L$. First, from \ref{bilingen} and the definitions of $\Delta^G(k;\textbf{i})(a,\textbf{Y})$, $\Delta^L(k;\textbf{i})(\textbf{Y})$, we obtain the following proposition, which implies that the minor $\Delta^G(k;\textbf{i})(a,\textbf{Y})$ is 
immediately obtained from $\Delta^L(k;\textbf{i})(\textbf{Y})$:
\begin{prop}\label{gprop}
We set $d:=i_k$ and $a=t^{h}$ $(h\in\mathfrak{h},\ t\in\mathbb{C}^{\times})$. Then
\[ 
\Delta^G(k;{\rm \bf{i}})(a,\textbf{Y})=
t^{(u_{\leq k}\Lambda_d)(h)}\Delta^L(k;{\rm \bf{i}})(\textbf{Y}).
\] 
\end{prop}
Next lemma will be used to lead the explicit formulas of $\Delta^G$ and $\Delta^L$:

\begin{lem}\label{gllem}
We suppose that $1\leq d<r$. If $m'+d>r$, then we have
\[ \Delta^G(k;{\rm \bf{i}})(a,\textbf{Y})=\lan  ax^L_{{\rm \bf{i}}}(\textbf{Y})(v_1\wedge v_2\wedge\cdots \wedge v_d) ,\q v_{m'+1}\wedge\cdots\wedge v_r\wedge v_{\ovl{d-r+m'}}\wedge\cdots\wedge v_{\ovl{1}}  \ran, \]
\[ \Delta^L(k;{\rm \bf{i}})(\textbf{Y})=\lan  x^L_{{\rm \bf{i}}}(\textbf{Y})(v_1\wedge v_2\wedge\cdots \wedge v_d) ,\q v_{m'+1}\wedge\cdots\wedge v_r\wedge v_{\ovl{d-r+m'}}\wedge\cdots\wedge v_{\ovl{1}}  \ran, \]
where $\lan ,\ran$ is the bilinear form we defined in \ref{bilingen}. If $m'+d\leq r$, we have
\begin{equation}\label{bglabel}
 \Delta^G(k;{\rm \bf{i}})(a,\textbf{Y})=\lan  ax^L_{{\rm \bf{i}}}(\textbf{Y})(v_1\wedge v_2\wedge\cdots \wedge v_d) ,\q v_{m'+1}\wedge\cdots\wedge v_{m'+d}  \ran, 
\end{equation}
\begin{equation}\label{bllabel}
 \Delta^L(k;{\rm \bf{i}})(\textbf{Y})=\lan  x^L_{{\rm \bf{i}}}(\textbf{Y})(v_1\wedge v_2\wedge\cdots \wedge v_d) ,\q v_{m'+1}\wedge\cdots\wedge v_{m'+d}  \ran. 
\end{equation}
\end{lem}
\nd
{\sl Proof.}

Let us prove this lemma for $\Delta^L(k;\textbf{i})(\textbf{Y})$ since the
case for $\Delta^G(k;\textbf{i})(a,\textbf{Y})$ 
is proven similarly.
Using (\ref{minor-bilin}) and (\ref{inc}), we see that $\Delta^L(k;\textbf{i})(\textbf{Y})=\Delta_{u_{\leq k}\Lambda_d, \Lambda_d}(x^L_{\textbf{i}}(\textbf{Y}))$ is given as 
\begin{equation}\label{bilin2}
\lan x^L_{\textbf{i}}(\textbf{Y})(v_1\wedge v_2\wedge\cdots \wedge v_d) ,\q \underbrace{\ovl{s_1}\cdots \ovl{s_r}}_{1{\rm \,st\
 cycle}}\cdots\underbrace{\ovl{s_1}\cdots \ovl{s_d}}_{m'{\rm \,th\ cycle}} (v_1\wedge v_2\wedge\cdots \wedge v_d) \ran.
\end{equation}
By (\ref{B-f1}), (\ref{B-f2}) and (\ref{smpl}), for $1\leq i \leq r-1$ and $1\leq j \leq r$, we get
\[ \ovl{s_i} v_j=
\begin{cases}
v_{i+1} & {\rm if}\ j=i, \\
-v_{i} & {\rm if}\ j=i+1, \\
v_j & {\rm if}\ {\rm otherwise},
\end{cases}\q 
\ovl{s_i} v_{\ovl{j}}=
\begin{cases}
v_{\ovl{i}} & {\rm if}\ j=i+1, \\
-v_{\ovl{i+1}} & {\rm if}\ j=i, \\
v_{\ovl{j}} & {\rm if}\ {\rm otherwise},
\end{cases}
\]
and we obtain
\[ \ovl{s_r} v_j=
\begin{cases}
v_{\ovl{r}} & {\rm if}\ j=r, \\
v_j & {\rm if}\ 1\leq j< r,
\end{cases}\q 
\ovl{s_r} v_{\ovl{j}}=
\begin{cases}
-v_{r} & {\rm if}\ j=r, \\
v_{\ovl{j}} & {\rm if}\ 1\leq j< r.
\end{cases}
\]
Therefore, if $m'+d\leq r$, then
\begin{equation}\label{ukeq1}u_{\leq k}(v_1\wedge\cdots\wedge v_d)=
\underbrace{\ovl{s_1}\cdots \ovl{s_r}}_{1{\rm \,st\ 
 cycle}}\cdots\underbrace{\ovl{s_1}\cdots \ovl{s_d}}_{m'{\rm \,th\ cycle}} (v_1\wedge \cdots \wedge v_d)=v_{m'+1}\wedge v_{m'+2}\wedge\cdots \wedge v_{m'+d}. 
\end{equation}
If $m'+d> r$, then we get
\begin{eqnarray}\label{ukeq2}
& &u_{\leq k}(v_1\wedge\cdots\wedge v_d)\nonumber \\
&=&
\underbrace{\ovl{s_1}\cdots \ovl{s_r}}_{1{\rm \,st\
 cycle}}\cdots\underbrace{\ovl{s_1}\cdots \ovl{s_d}}_{m'{\rm \,th\ cycle}} (v_1\wedge v_2\wedge\cdots \wedge v_d) \nonumber \\
&=&\underbrace{\ovl{s_1}\cdots \ovl{s_r}}_{1{\rm \,st
\ cycle}}\cdots\underbrace{\ovl{s_1}\cdots \ovl{s_r}}_{m'-r+d{\rm \,th\ cycle}} (v_{r-d+1}\wedge\cdots \wedge v_{r}) \nonumber \\
&=&\underbrace{\ovl{s_1}\cdots \ovl{s_r}}_{1{\rm \,st
\ cycle}}\cdots\underbrace{\ovl{s_1}\cdots \ovl{s_r}}_{m'-r+d-1{\rm \,th\ cycle}} (v_{r-d+2}\wedge\cdots \wedge v_{r}\wedge v_{\ovl{1}})\nonumber \\
&=&\underbrace{\ovl{s_1}\cdots \ovl{s_r}}_{1{\rm \,st
\ cycle}}\cdots\underbrace{\ovl{s_1}\cdots \ovl{s_r}}_{m'-r+d-2{\rm \,th\ cycle}} (v_{r-d+3}\wedge\cdots \wedge v_{r}\wedge v_{\ovl{1}}\wedge (-v_{\ovl{2}}))\nonumber \\
&=&\cdots\q =v_{m'+1}\wedge\cdots \wedge v_{r}\wedge v_{\ovl{1}}\wedge (-v_{\ovl{2}})\wedge 
 ((-1)^2v_{\ovl{3}})\wedge\cdots\wedge ((-1)^{d-r+m'-1}v_{\ovl{d-r+m'}})
 \nonumber \\
&=&v_{m'+1}\wedge\cdots \wedge v_{r}\wedge v_{\ovl{d-r+m'}}\wedge\cdots\wedge v_{\ovl{1}}.
\end{eqnarray}
Hence, we get our claim by (\ref{bilin2}). \qed 

\vspace{3mm}

In the end of this subsection, let us prove the following lemma:

\begin{lem}\label{gmlem}
Let $u$, ${\rm \bf{i}}$ and $\textbf{Y}$ be as in $(\ref{uvset})$, $(\ref{iset})$ and $(\ref{yset})$. Let $i_{n+1}\in [1,r]$ be an index such that $u':=us_{i_{n+1}}\in W$ satisfies $l(u')>l(u)$. We set the reduced word ${\rm \bf{i}}'$ for $u'$ as
\[
{\rm \bf{i}}'=(\underbrace{1,\cdots,r}_{1{\rm \,st\
  cycle}},\underbrace{1,\cdots,r}_{2{\rm \,nd\
  cycle}},\cdots,\underbrace{1,\cdots,r}_{m-1{\rm \,th\ cycle}},
\underbrace{1,\cdots,i_n}_{m{\rm \,th\ cycle}},i_{n+1}),
\]
and write $\textbf{Y}'\in(\mathbb{C}^{\times})^{n+1}$ as
\[ \textbf{Y}':=(Y_{1,1},\cdots,Y_{1,r},\cdots,Y_{m-1,1},\cdots,Y_{m-1,r},Y_{m,1},\cdots,Y_{m,i_n},Y).
\] 
For an integer $k$ $(1\leq k\leq n)$, if $d:=i_k\neq i_{n+1}$, then $\Delta^L(k;{\rm \bf{i}}')(\textbf{Y}')$ does not depend on $Y$, so we can regard it as a function on $(\mathbb{C}^{\times})^{n}$. Furthermore, we have 
\begin{equation}\label{Lomit}
\Delta^L(k;{\rm \bf{i}})(\textbf{Y})=\Delta^L(k;{\rm \bf{i}}')(\textbf{Y}') .
\end{equation}
\end{lem}
\nd
{\sl Proof.}
By the definition (\ref{xldef}) of $x^L_{\textbf{i}}$, we have
\begin{equation}\label{gmlempr1}
x^L_{\textbf{i}'}(\textbf{Y}')=x^L_{\textbf{i}}(\textbf{Y})x_{-i_{n+1}}(Y). 
\end{equation}
It follows from $x_{-i_{n+1}}(Y):={\rm exp}(Y f_{i_{n+1}})\cdot(Y^{-h_{i_{n+1}}})$ (see (\ref{alxmdef})) that for $j$ $(1\leq j\leq r-1)$,
\begin{equation}\label{xmpro}
x_{-i_{n+1}}(Y)v_j=
\begin{cases}
Y^{-1}v_{i_{n+1}}+v_{i_{n+1}+1} & {\rm if}\ j=i_{n+1}, \\
Y v_{i_{n+1}+1} & {\rm if}\ j=i_{n+1}+1, \\
v_j & {\rm otherwise}.
\end{cases}
\end{equation}
Thus, if $d<i_{n+1}$, then we have $x_{-i_{n+1}}(Y)(v_1\wedge\cdots\wedge v_d)=v_1\wedge\cdots\wedge v_d$.
If $d>i_{n+1}$, then we have
\begin{flushleft}
$x_{-i_{n+1}}(Y)(v_1\wedge\cdots\wedge v_d)$
\end{flushleft}
\begin{eqnarray*}
&=&v_1\wedge\cdots \wedge v_{i_{n+1}-1}\wedge (Y^{-1}v_{i_{n+1}}+v_{i_{n+1}+1})\wedge Y v_{i_{n+1}+1} \wedge\cdots \wedge v_d \\
&=&v_1\wedge\cdots\wedge v_d.
\end{eqnarray*}
Since we assume $i_{n+1}\neq d$, we get
\begin{equation}\label{gmlempr2}
x_{-i_{n+1}}(Y)(v_1\wedge\cdots\wedge v_d)=v_1\wedge\cdots\wedge v_d.
\end{equation}
We can easily see that $u_{\leq k}=u'_{\leq k}(=\underbrace{s_1\cdots s_r}_{1{\rm \,st\
 cycle}}\cdots\underbrace{s_1\cdots s_d}_{m'{\rm \,th\ cycle}})$. Therefore, it follows from (\ref{minor-bilin}), (\ref{gmlempr1}) and (\ref{gmlempr2}) that
\begin{eqnarray*}
\Delta^L(k;\textbf{i}')(\textbf{Y}')&=&\Delta_{u'_{\leq k}\Lambda_d,\Lambda_d}(x^L_{\textbf{i}'}(\textbf{Y}')) \\
&=&\lan x^L_{\textbf{i}'}(\textbf{Y}')(v_1\wedge v_2\wedge\cdots \wedge v_d) ,\ u'_{\leq k} (v_1\wedge v_2\wedge\cdots \wedge v_d) \ran \\
&=&\lan x^L_{\textbf{i}}(\textbf{Y})(v_1\wedge \cdots \wedge v_d) ,\ u_{\leq k} (v_1\wedge \cdots \wedge v_d) \ran=\Delta^L(k;\textbf{i})(\textbf{Y}),
\end{eqnarray*}
which is our desired result. \qed

\vspace{5mm}

In the rest of the paper, we will only treat $\Delta^L(k;\textbf{i})(\textbf{Y})$ and assume that $i_n=i_k$
due to Proposition \ref{gprop} and this lemma. 

\subsection{Generalized minor $\Delta^L(k;\textbf{i})(\textbf{Y})$ of type ${\rm B}_r$}\label{Mainsec}

Now, we shall give the explicit formula of generalized minors $\Delta^L(k;\textbf{i})(\textbf{Y})$ on reduced double Bruhat cells $L^{u,e}$ of type ${\rm B}_r$. For $l\in[1,m]$ and $1\leq k\leq \ovl{1}$, we set the Laurent monomials
\begin{equation}\label{bbbar}
B(l,k):=
\begin{cases}
\frac{Y_{l,k-1}}{Y_{l,k}} & {\rm if}\ 1\leq k\leq r-1, \\
\frac{Y_{l,r-1}}{Y^2_{l,r}} & {\rm if}\ k= r, \\
\frac{Y_{l,r}}{Y_{l+1,r}} & {\rm if}\ k= 0, \\
\frac{Y^2_{l,r}}{Y_{l+1,r-1}} & {\rm if}\ k= \ovl{r}, \\
\frac{Y_{l,|k|}}{Y_{l+1,|k|-1}} & {\rm if}\ \ovl{r-1}\leq k\leq \ovl{1}.
\end{cases}
\end{equation}

For $1\leq l\leq r$, we set $|l|=|\ovl{l}|=l$. Then the following theorem holds:

\begin{thm}\label{thm1}
In the above setting, we set
\[ \textbf{Y}:=(Y_{1,1},Y_{1,2},\cdots,Y_{1,r},\cdots,Y_{m-1,1},\cdots,Y_{m-1,r},Y_{m,1},\cdots,Y_{m,i_n})\in  (\mathbb{C}^{\times})^{n}, \]
$d:=i_k=i_n$ and suppose that $d<r$. Then we have
\begin{multline*} 
\Delta^L(k;{\rm \bf{i}})(\textbf{Y})= \sum_{(*)} 2^{C[\{k^{(s)}_i\}]}\prod^{d}_{i=1} B(m-l^{(1)}_i,k^{(1)}_i)B(m-l^{(2)}_i,k^{(2)}_i)\cdots \\
\cdots  B(m-l^{(m-m'-1)}_i,k^{(m-m'-1)}_i) B(m-l^{(m-m')}_i,k^{(m-m')}_i),
\end{multline*} 
where
\[ l^{(s)}_{i}:=
\begin{cases}
k^{(s)}_{i}+s-i-1 & {\rm if}\ k^{(s)}_{i}\in\{1,2,\cdots,r\}, \\
s-i+r & {\rm if}\ k^{(s)}_{i}\in \{0,\ovl{r},\ovl{r-1},\cdots,\ovl{1}\},
\end{cases}
\q (1\leq i\leq d)
\]
\begin{multline*}
C[\{k^{(s)}_i\}]=C[\{k^{(s)}_i\}_{1\leq i\leq d,\ 1\leq s\leq m-m'}]:=\\
\#\{j\in[1,d-1]|\ k^{(t)}_{j+1}=0,\ k^{(t+1)}_{j+1}\neq 0\ {\rm and}\ k^{(t)}_j\neq0\ {\rm for\ some}\ t\in[1,m-m']\},
\end{multline*}
and $(*)$ is the conditions for $k^{(s)}_i$ $(1\leq s\leq m-m',\ 1\leq i\leq d)$ $:$
\begin{equation}\label{B-cond1}
1\leq k^{(s)}_1\leq k^{(s)}_2\leq \cdots\leq k^{(s)}_d\leq\ovl{1},\q k^{(s)}_i=k^{(s)}_{i+1}\ {\rm if\ and\ only\ if}\ k^{(s)}_i=k^{(s)}_{i+1}=0,
\end{equation}
\begin{equation}\label{B-cond2}
1\leq k^{(1)}_i\leq \cdots\leq k^{(m-m')}_i\leq m'+i\q (1\leq i\leq r-m'),
\end{equation}
\begin{equation}\label{B-cond3}
1\leq k^{(1)}_i\leq \cdots\leq k^{(m-m')}_i\leq \ovl{1}\q (r-m'+1\leq i\leq d),
\end{equation}
\begin{equation}\label{B-cond4}
k^{(s)}_i<k^{(s+1)}_{i+1}.
\end{equation}
\end{thm}

Next, we set 
\[ B(l,r+1):=\frac{1}{Y_{l,r}}.\]
\begin{thm}\label{thm2}
In the above setting, we set
\[ \textbf{Y}:=(Y_{1,1},Y_{1,2},\cdots,Y_{1,r},\cdots,Y_{m-1,1},\cdots,Y_{m-1,r},Y_{m,1},\cdots,Y_{m,i_n})\in  (\mathbb{C}^{\times})^{n}, \]
$d:=i_k=i_n$ and suppose that $d=r$. Then we have
\begin{eqnarray*} \Delta^L(k;{\rm \bf{i}})(\textbf{Y})&=&\sum_{(*)}\prod^{m'+1}_{i=1} B(t_i-1,\ovl{k^{(t_i-1)}_i})B(t_i-2,\ovl{k^{(t_i-2)}_i})\cdots B(m'-i+1,\ovl{k^{(m'-i+1)}_i})\\ 
& &\qq \qq \cdot
B(t_{i-1}-1,r+1)B(t_{i-1}-2,r+1)\cdots B(t_{i}+1,r+1),
\end{eqnarray*}
where $(*)$ is the following conditions:
\[ 1\leq t_{m'}<t_{m'-1}<\cdots<t_1\leq m,\ t_0=m+1,\ t_{m'+1}=0, \]
\[ i\leq k^{(m'-i+1)}_i\leq\cdots \leq k^{(t_i-2)}_i\leq k^{(t_i-1)}_i\leq r, \]
\[ k^{(s-1)}_i\leq k^{(s)}_i<k^{(s-1)}_{i+1},\ \ (1\leq s\leq t_i-1). \]
\end{thm}

\begin{ex}\label{pathex3}
For {\rm rank} $r=3$, $u=s_1s_2s_3s_1s_2s_3s_1s_2$, $k=5$ 
and the reduced word ${\rm \bf{i}}=(1,2,3,1,2,3,1,2)$ for $u$, 
we have $m=3$, $m'=2$ and $d=2$ $($see $(\ref{uvset}),\ (\ref{iset}))$. 
Then, we have $s=1$ and we write $k_i$ for $k^{(s)}_i$. 
Thus, the set of  all $(k_1,k_2)$ satisfying $(*)$ in Theorem \ref{thm1} is
\begin{multline*}
\{(1,2),(1,3),(1,0),(1,\ovl 3),(1,\ovl 2),(1,\ovl 1),(2,3),(2,0),(2,\ovl 3),(2,\ovl 2), \\
(2,\ovl 1),(3,0),(3,\ovl 3),(3,\ovl 2),(3,\ovl 1)\}
\end{multline*}
Here, for all $(k_1,k_2)$ the corresponding monomials are as follows:
\[
\begin{array}{ccc}
(1,2)\leftrightarrow B(3,1)B(3,2)&
(1,3)\leftrightarrow    B(3,1)B(2,3)&
(1,0)\leftrightarrow B(3,1)B(1,0)\\
(1,\ovl 3)\leftrightarrow B(3,1)B(1,\ovl{3}) &
(1,\ovl 2)\leftrightarrow B(3,1)B(1,\ovl{2})&
(1,\ovl 1)\leftrightarrow  B(3,1)B(1,\ovl{1})\\
(2,3)\leftrightarrow  B(2,2)B(2,3)&
(2,0)\leftrightarrow B(2,2)B(1,0) &
(2,\ovl 3)\leftrightarrow B(2,2)B(1,\ovl{3})\\
(2,\ovl 2)\leftrightarrow B(2,2)B(1,\ovl{2})&
(2,\ovl 1)\leftrightarrow B(2,2)B(1,\ovl{1})&
(3,0)\leftrightarrow B(1,3)B(1,0)\\
(3,\ovl 3)\leftrightarrow B(1,3)B(1,\ovl{3})&
(3,\ovl 2)\leftrightarrow B(1,3)B(1,\ovl{2})&
(3,\ovl 1)\leftrightarrow B(1,3)B(1,\ovl{1})\\
\end{array}
\]
Thus, we obtain:
\begin{eqnarray*}
\Delta^L(5;{\rm \bf{i}})(\textbf{Y})&=&
B(3,1)B(3,2)
+B(3,1)B(2,3)
+2B(3,1)B(1,0)
+B(3,1)B(1,\ovl{3}) \\
& &
+B(3,1)B(1,\ovl{2})
+B(3,1)B(1,\ovl{1})
+B(2,2)B(2,3)
+2B(2,2)B(1,0) \\
& &
+B(2,2)B(1,\ovl{3})
+B(2,2)B(1,\ovl{2})
+B(2,2)B(1,\ovl{1})
+2B(1,3)B(1,0)  \\
& &
+B(1,3)B(1,\ovl{3})+B(1,3)B(1,\ovl{2})+B(1,3)B(1,\ovl{1})\\
&=&
\frac{1}{Y_{3,2}}+\frac{Y_{2,2}}{Y_{3,1}Y^2_{2,3}}+2\frac{Y_{1,3}}{Y_{3,1}Y_{2,3}}+\frac{Y^2_{1,3}}{Y_{3,1}Y_{2,2}}
+\frac{Y_{1,2}}{Y_{3,1}Y_{2,1}}+\frac{Y_{1,1}}{Y_{3,1}}+\frac{Y_{2,1}}{Y^2_{2,3}}\\
\q &+&2\frac{Y_{2,1}Y_{1,3}}{Y_{2,3}Y_{2,2}}+\frac{Y_{2,1}Y^2_{1,3}}{Y^2_{2,2}}+2\frac{Y_{1,2}}{Y_{2,2}}+
\frac{Y_{2,1}Y_{1,1}}{Y_{2,2}}+2\frac{Y_{1,2}}{Y_{2,3}Y_{1,3}}+\frac{Y^2_{1,2}}{Y^2_{1,3}Y_{2,1}}+\frac{Y_{1,2}Y_{1,1}}{Y^2_{1,3}}.
\end{eqnarray*}
\end{ex}

\section{The proof of Theorem \ref{thm1}\ and \ref{thm2}}\label{prsec}

In this section, we shall give the proof of Theorem \ref{thm1} and \ref{thm2}.

\subsection{The set $X_d(m,m')$ of paths}

In this subsection, we shall introduce a set $X_d(m,m')$ of ``paths''
which will correspond to the set of the terms in $\Delta^L(k;\textbf{i})(\textbf{Y})$. Let
$m$, $m'$ and $d$ $(d<r)$ be the positive integers as in 
\ref{Mainsec}. We set $J:=\{j,\ovl{j}|\ 1\leq j\leq r\}\cup\{0\}$ and for $l\in\{1,2,\cdots,r\}$, set $|l|=|\ovl{l}|=l$.

\begin{defn}\label{vdeddef}
Let us define the directed graph $(V_d,E_d)$ as follows:
The set $V_d=V_d(m)$ of vertices is defined by 
\[V_d(m):=\{{\rm vt}(m-s;a^{(s)})|\ a^{(s)}=(a^{(s)}_1,a^{(s)}_2,\cdots,a^{(s)}_d)\in J^d,\ 0\leq s\leq m \}. 
\]
And we define the set $E_d=E_d(m)$ of directed edges as 
\begin{multline*} E_d(m):=\{{\rm vt}(m-s;a^{(s)})\rightarrow
{\rm vt}(m-s-1;a^{(s+1)})\\
|\ {\rm vt}(m-s;a^{(s)}),\ {\rm vt}(m-s-1;a^{(s+1)})\in V_d(m),\ 0\leq s\leq m-1\}.
\end{multline*}
\end{defn}

Now, let us define the set of directed paths from ${\rm vt}(m;1,2,\cdots,d)$ to ${\rm vt}(0;m'+1,m'+2,\cdots,r,\ovl{d-r+m'},\ovl{d-r+m'-1},\cdots,\ovl{2},\ovl{1})$ (resp. ${\rm vt}(0;m'+1,m'+2,\cdots,m'+d)$) in the case $m'+d>r$ (resp. $m'+d\leq r$) in $(V_d,E_d)$.

\begin{defn}\label{pathdef}
Let $X_d(m,m')$ be the set of directed paths $p$
\begin{multline*}p={\rm vt}(m;a^{(0)}_1,\cdots,a^{(0)}_d)
\rightarrow{\rm vt}(m-1;a^{(1)}_1,\cdots,a^{(1)}_d)\rightarrow {\rm vt}(m-2;a^{(2)}_1,
\cdots,a^{(2)}_d)\\
\rightarrow
\cdots\rightarrow{\rm vt}(1;a^{(m-1)}_1,\cdots,a^{(m-1)}_d)
\rightarrow{\rm vt}(0;a^{(m)}_1,\cdots,a^{(m)}_d),
\end{multline*}
which satisfy the following conditions: For $s\in\mathbb{Z}$ $(0\leq s\leq m)$,
\begin{enumerate}
\item $a^{(s)}_{\zeta}\in J$ $(1\leq \zeta\leq d)$, 
\item $a^{(s)}_{1}<a^{(s)}_{2}<\cdots<a^{(s)}_{d}$ in the order (\ref{B-order}), 
\item If $a^{(s)}_{\zeta}\in\{j|1\leq j\leq r-1\}$, then $a^{(s+1)}_{\zeta}=a^{(s)}_{\zeta}$ or $a^{(s)}_{\zeta}+1$. If $a^{(s)}_{\zeta}\in\{r,\ovl{r},\ovl{r-1},\cdots,\ovl{1}\}\cup\{0\}$, then $a^{(s)}_{\zeta}\leq a^{(s+1)}_{\zeta}\leq\ovl{1}$ in the order (\ref{B-order}),
\item $(a^{(0)}_1,a^{(0)}_2,\cdots,a^{(0)}_d)=(1,2,\cdots,d)$ and
\[ (a^{(m)}_1,\cdots,a^{(m)}_d)=
\begin{cases}
(m'+1,m'+2,\cdots,r,\ovl{d-r+m'},\cdots,\ovl{2},\ovl{1}) & {\rm if}\ m'+d>r,\\
(m'+1,m'+2,\cdots,m'+d) & {\rm if}\ m'+d\leq r,
\end{cases}
\]
\item If $a^{(s+1)}_{\zeta}\in\{\ovl{r},\ovl{r-1},\cdots,\ovl{1}\}\cup\{0\}$, then $a^{(s+1)}_{\zeta}\leq a^{(s)}_{\zeta+1}$ in the order (\ref{B-order}). Furthermore, $a^{(s+1)}_{\zeta}= a^{(s)}_{\zeta+1}$ if and only if $a^{(s+1)}_{\zeta}= a^{(s)}_{\zeta+1}=0$.
\end{enumerate}
\end{defn}

\begin{defn}\label{connected-def}
We say that two vertices ${\rm vt}(m-s;a^{(s)}_1,\cdots,a^{(s)}_d)$ and ${\rm vt}(m-s-1;a^{(s+1)}_1,\cdots,a^{(s+1)}_d)$ are {\it connected} if these vertices satisfy the conditions (i), (ii), (iii) and (v) in Definition \ref{pathdef}.
\end{defn}


Define a Laurent monomial associated with each edge of paths in $X_d(m,m')$.
\begin{defn}\label{labeldef}
\begin{enumerate}
\item For each $s$ $(0\leq s\leq m-1)$ and $i\in[1,r-1]$, we set
\[Q^{(s)}(i\rightarrow j):=
\begin{cases}
\frac{Y_{m-s,i-1}}{Y_{m-s,i}} & {\rm if}\ j=i,\\
1 & {\rm if}\ j=i+1,
\end{cases} 
\]
\[Q^{(s)}(r\rightarrow j):=
\begin{cases}
\frac{Y_{m-s,r-1}}{Y^2_{m-s,r}} & {\rm if}\ j=r,\\
\frac{1}{Y_{m-s,r}} & {\rm if}\ j=0, \\
\frac{1}{Y_{m-s,|j|-1}} & {\rm if}\ j\in\{\ovl{r},\ovl{r-1}\cdots\ovl{1}\}, \\
\end{cases} 
\]
and
\[Q^{(s)}(0\rightarrow j):=
\begin{cases}
1 & {\rm if}\ j=0, \\
2\frac{Y_{m-s,r}}{Y_{m-s,|j|-1}} & {\rm if}\ j\in\{\ovl{r},\ovl{r-1}\cdots\ovl{1}\}. \\
\end{cases} 
\]
For $i,j\in\{\ovl{r},\ovl{r-1}\cdots\ovl{1}\}$ $(i\leq j)$, we also set
\[Q^{(s)}(\ovl{r}\rightarrow j):=\frac{Y^2_{m-s,r}}{Y_{m-s,|j|-1}},\ \ Q^{(s)}(i\rightarrow j):=\frac{Y_{m-s,|i|}}{Y_{m-s,|j|-1}}.\]
For an edge $e^{i_1,\cdots,i_d}_{j_1,\cdots,j_d}={\rm vt}(m-s;i_1,\cdots,i_d)\rightarrow {\rm vt}(m-s-1;j_1,\cdots,j_d)$ in $E_d(m)$, we define the {\it label} $Q^{(s)}(e^{i_1,\cdots,i_d}_{j_1,\cdots,j_d})$ of the edge $e^{i_1,\cdots,i_d}_{j_1,\cdots,j_d}$ as
\[ Q^{(s)}(e^{i_1,\cdots,i_d}_{j_1,\cdots,j_d}):=c^{(s)}\prod^d_{k=1}Q^{(s)}(i_k\rightarrow j_k), \]
where
\[
c^{(s)}:=
\begin{cases}
\frac{1}{2} & {\rm if\ there\ exists\ some\ } k\ {\rm such\ that}\ j_k=i_{k+1}=0, \\
1 & {\rm otherwise}. 
\end{cases}
\]

\item Let $p\in X_d(m,m')$ be a path:
\begin{multline*}p={\rm vt}(m;a^{(0)}_1,\cdots,a^{(0)}_d)
\rightarrow{\rm vt}(m-1;a^{(1)}_1,\cdots,a^{(1)}_d)\rightarrow {\rm vt}(m-2;a^{(2)}_1,
\cdots,a^{(2)}_d)\\
\rightarrow
\cdots\rightarrow{\rm vt}(1;a^{(m-1)}_1,\cdots,a^{(m-1)}_d)
\rightarrow{\rm vt}(0;a^{(m)}_1,\cdots,a^{(m)}_d).
\end{multline*}
For each $s$ $(0\leq s\leq m-1)$, we denote the label of the $(m-s)$ th edge ${\rm vt}(m-s;a^{(s)}_1,a^{(s)}_2,\cdots$, $a^{(s)}_d)\rightarrow{\rm vt}(m-s-1;a^{(s+1)}_1,a^{(s+1)}_2,\cdots,a^{(s+1)}_d)$ of $p$ by $Q^{(s)}(p)$. And we define the {\it label} $Q(p)$ {\it of the path} $p$ as the total product:
\begin{equation}\label{alllabels} Q(p):=\prod_{s=0}^{m-1}Q^{(s)}(p). 
\end{equation}
\item For $s,s'\in\mathbb{Z}$ $(0\leq s'<s''\leq m)$, if the vertices of the path $p'$
\[
p'={\rm vt}(m-s';a^{(s')})\rightarrow{\rm vt}(m-s'-1;a^{(s'+1)})\rightarrow 
\cdots
\rightarrow{\rm vt}(m-s'';a^{(s'')})
\]
satisfy the (i),(ii),(iii) and (v) in Definition \ref{pathdef}, we call $p'$ {\it subpath} from ${\rm vt}(m-s';a^{(s')})$ to ${\rm vt}(m-s'';a^{(s'')})$, and define the {\it label} {\it of} $p'$ as 
\begin{equation}\label{sublabels} 
Q(p'):=\prod_{s=s'}^{s''-1}Q^{(s)}(p'),
\end{equation}
where $Q^{(s)}(p')$ is the label of the edge ${\rm vt}(m-s;a^{(s)})\rightarrow{\rm vt}(m-s-1;a^{(s+1)})$.
\end{enumerate}
\end{defn}

\begin{ex}\label{pathex}
Let $r=m=3$, $m'=2$, $d=2$. We can describe the paths of $X_2(3,2)$ as follows. For simplicity, we denote vertices ${\rm vt}(s;a_1,a_2)$ by $(a_1,a_2)$. For example, the top vertex $(1,2)$ implies ${\rm vt}(3;1,2)$ and second vertices $(1,2)$, $(1,3)$ and $(2,3)$ imply ${\rm vt}(1,2)$, ${\rm vt}(1,3)$ and ${\rm vt}(2,3)$ respectively.

\begin{figure}[htbp]\label{Figure1}
\begin{xy}
(0,90) *{(1,2)}="3;1,2",
(0,70)*{(1,2)}="2;1,2",
(15,70)*{(1,3)}="2;1,3",
(30,70)*{(2,3)}="2;2,3",
(30,50)*{(2,3)}="1;2,3",
(42,50)*{(2,0)}="1;2,0",
(54,50)*{(2,\ovl{3})}="1;2,4",
(66,50)*{(2,\ovl{2})}="1;2,5",
(78,50)*{(2,\ovl{1})}="1;2,6",
(90,50)*{(3,0)}="1;3,0",
(102,50)*{(3,\ovl{3})}="1;3,4",
(114,50)*{(3,\ovl{2})}="1;3,5",
(126,50)*{(3,\ovl{1})}="1;3,6",
(126,30)*{(3,\ovl{1})}="0;3,4",
\ar@{->} "3;1,2";"2;1,2"
\ar@{->} "3;1,2";"2;1,3"
\ar@{->} "3;1,2";"2;2,3"
\ar@{->} "2;1,2";"1;2,3"
\ar@{->} "2;1,3";"1;2,3"
\ar@{->} "2;1,3";"1;2,4"
\ar@{->} "2;1,3";"1;2,0"
\ar@{->} "2;1,3";"1;2,5"
\ar@{->} "2;1,3";"1;2,6"
\ar@{->} "2;2,3";"1;2,3"
\ar@{->} "2;2,3";"1;2,0"
\ar@{->} "2;2,3";"1;2,4"
\ar@{->} "2;2,3";"1;2,5"
\ar@{->} "2;2,3";"1;2,6"
\ar@{->} "2;2,3";"1;3,0"
\ar@{->} "2;2,3";"1;3,4"
\ar@{->} "2;2,3";"1;3,5"
\ar@{->} "2;2,3";"1;3,6"
\ar@{->} "1;2,3";"0;3,4"
\ar@{->} "1;2,0";"0;3,4"
\ar@{->} "1;2,4";"0;3,4"
\ar@{->} "1;2,5";"0;3,4"
\ar@{->} "1;2,6";"0;3,4"
\ar@{->} "1;3,0";"0;3,4"
\ar@{->} "1;3,4";"0;3,4"
\ar@{->} "1;3,5";"0;3,4"
\ar@{->} "1;3,6";"0;3,4"
\end{xy}\caption{The paths in $X_2(3,2)$}
\end{figure}
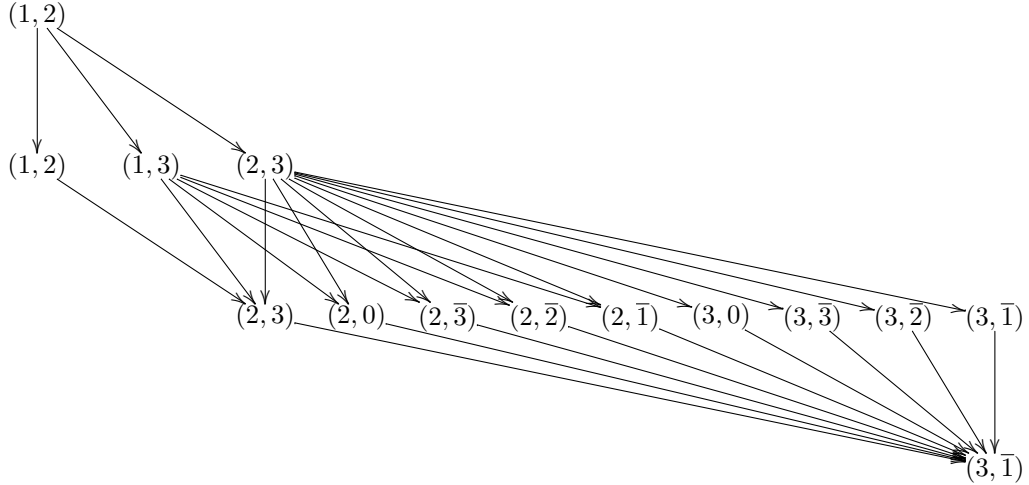
Thus, $X_2(3,2)$ has the following paths$:$

$p_1=(1,2)\rightarrow(1,2)\rightarrow(2,3)\rightarrow(3,\ovl{1})$,\ \  $p_2=(1,2)\rightarrow(1,3)\rightarrow(2,3)\rightarrow(3,\ovl{1})$,

$p_3=(1,2)\rightarrow(1,3)\rightarrow(2,0)\rightarrow(3,\ovl{1})$,\ \ 
$p_4=(1,2)\rightarrow(1,3)\rightarrow(2,\ovl{3})\rightarrow(3,\ovl{1})$,

$p_5=(1,2)\rightarrow(1,3)\rightarrow(2,\ovl{2})\rightarrow(3,\ovl{1})$,\ \ 
$p_6=(1,2)\rightarrow(1,3)\rightarrow(2,\ovl{1})\rightarrow(3,\ovl{1})$,
 
$p_7=(1,2)\rightarrow(2,3)\rightarrow(2,3)\rightarrow(3,\ovl{1})$,\ \ 
$p_8=(1,2)\rightarrow(2,3)\rightarrow(2,0)\rightarrow(3,\ovl{1})$,

$p_9=(1,2)\rightarrow(2,3)\rightarrow(2,\ovl{3})\rightarrow(3,\ovl{1})$,\ 
$p_{10}=(1,2)\rightarrow(2,3)\rightarrow(2,\ovl{2})\rightarrow(3,\ovl{1})$,

$p_{11}=(1,2)\rightarrow(2,3)\rightarrow(2,\ovl{1})\rightarrow(3,\ovl{1})$,

$p_{12}=(1,2)\rightarrow(2,3)\rightarrow(3,0)\rightarrow(3,\ovl{1})$,

$p_{13}=(1,2)\rightarrow(2,3)\rightarrow(3,\ovl{3})\rightarrow(3,\ovl{1})$,

$p_{14}=(1,2)\rightarrow(2,3)\rightarrow(3,\ovl{2})\rightarrow(3,\ovl{1})$,

$p_{15}=(1,2)\rightarrow(2,3)\rightarrow(3,\ovl{1})\rightarrow(3,\ovl{1})$.

\vspace{2mm}

Let us calculate the label of the path $p_1$. 
By Definition \ref{labeldef} $(ii)$, the label $Q^{(0)}(p_1)$ of the edge ${\rm vt}(3;1,2)\rightarrow{\rm vt}(2;1,2)$ is
\[ Q^{(0)}(p_1)=\frac{Y_{3,1-1}}{Y_{3,1}}\frac{Y_{3,2-1}}{Y_{3,2}}=\frac{1}{Y_{3,2}}, \] 
where we set $Y_{3,0}=1$ following Remark \ref{importantrem}. The labels of the edges ${\rm vt}(2;1,2)\rightarrow{\rm vt}(1;2,3)$ and ${\rm vt}(1;2,3)\rightarrow{\rm vt}(1;3,\ovl{1})$ are $Q^{(1)}(p_1)=1$,\ 
$Q^{(2)}(p_1)=1$. Therefore, we get $Q(p_1)=\frac{1}{Y_{3,2}}$. Similarly, we have
\[ Q(p_1)=\frac{1}{Y_{3,2}},\ \ Q(p_2)=\frac{Y_{2,2}}{Y_{3,1}Y^2_{2,3}},\ 
\ Q(p_3)=\frac{2Y_{1,3}}{Y_{3,1}Y_{2,3}},\ \ Q(p_4)=\frac{Y^2_{1,3}}{Y_{3,1}Y_{2,2}}, \]
\[ Q(p_5)=\frac{Y_{1,2}}{Y_{3,1}Y_{2,1}},\ \ Q(p_6)=\frac{Y_{1,1}}{Y_{3,1}},
\ \ Q(p_7)=\frac{Y_{2,1}}{Y_{2,3}^2},\ \ Q(p_8)=\frac{2Y_{2,1}Y_{1,3}}{Y_{2,2}Y_{2,3}}, \]
\[ Q(p_9)=\frac{Y_{2,1}Y^2_{1,3}}{Y^2_{2,2}},\ \ Q(p_{10})=\frac{Y_{1,2}}{Y_{2,2}},
\ \ Q(p_{11})=\frac{Y_{2,1}Y_{1,1}}{Y_{2,2}},\ \ Q(p_{12})=\frac{2Y_{1,2}}{Y_{2,3}Y_{1,3}}, \]
\[ Q(p_{13})=\frac{Y_{1,2}}{Y_{2,2}},\ \ Q(p_{14})=\frac{Y^2_{1,2}}{Y_{2,1}Y^2_{1,3}},\ \ Q(p_{15})=\frac{Y_{1,1}Y_{1,2}}{Y^2_{1,3}}. \]
\end{ex}

\begin{defn}\label{iseq}
For each path $p\in X_d(m,m')$
\begin{multline*}p={\rm vt}(m;a^{(0)}_1,\cdots,a^{(0)}_d)\rightarrow{\rm vt}(m-1;a^{(1)}_1,\cdots,a^{(1)}_d)\rightarrow{\rm vt}(m-2;a^{(2)}_1,\cdots,a^{(2)}_d)\\
\rightarrow
\cdots\rightarrow{\rm vt}(1;a^{(m-1)}_1,\cdots,a^{(m-1)}_d)\rightarrow{\rm vt}(0;a^{(m)}_1,\cdots,a^{(m)}_d)
\end{multline*}
and $i\in \{1,\cdots,d\}$, we call the sequence $a^{(0)}_i\rightarrow a^{(1)}_i\rightarrow a^{(2)}_i\rightarrow\cdots\rightarrow a^{(m)}_i$ an {\it i-sequence} of $p$.
\end{defn}

We can easily see the following by Definition \ref{pathdef} (iii) and (iv): For $i$ $(1\leq i\leq d)$,
\begin{equation}\label{iseqine}
\begin{cases}
i=a^{(0)}_i\leq a^{(1)}_i\leq \cdots\leq a^{(m)}_i\leq m'+i, & {\rm if}\ i\leq r-m', \\
i=a^{(0)}_i\leq a^{(1)}_i\leq \cdots\leq a^{(m)}_i\leq \ovl{d-i+1}, & {\rm if}\ r-m'<i\leq d, \\
\end{cases}
\end{equation}
in the order (\ref{B-order}).

\subsection{One-to-one correspondence between paths in $X_d(m,m')$ and terms of $\Delta^L(k;\textbf{i})(\textbf{Y})$}

In this section, we describe the terms in $\Delta^L(k;\textbf{i})(\textbf{Y})$ as the paths in $X_d(m,m')$:

\begin{prop}\label{pathlem}We use the setting and the notations in Sect.\ref{gmc}:
\[ u=(s_1s_2\cdots s_r)^{m-1}s_1\cdots s_{i_n},\q v=e.\]
Then, we have the following:
\[ \Delta^L(k;{\rm \bf{i}})(\textbf{Y})=\sum_{p\in X_d(m,m')} Q(p). \]
\end{prop}

Let us give an overview of the proof of Proposition \ref{pathlem}. For $1\leq s\leq m$, we define
\begin{equation}\label{xmdef1}
x^{(s)}_{-[1,r]}:=x_{-1}(Y_{s,1})\cdots x_{-r}(Y_{s,r}).
\end{equation}
For $1\leq s\leq m$ and $i_1,\cdots,i_d\in J=\{i,\ovl i|1\leq i\leq r\}\cup\{0\}$, we set 
\begin{equation}\label{xmdef2}
 (s;i_1,i_2,\cdots,i_d):=\lan x^{(1)}_{-[1,r]} x^{(2)}_{-[1,r]}\cdots x^{(s)}_{-[1,r]}
(v_{i_1}\wedge\cdots\wedge v_{i_d}), \ u_{\leq k}(v_1\wedge\cdots\wedge v_d)  \ran . 
\end{equation}

We shall prove $\Delta^L(k;\textbf{i})(\textbf{Y})=(m;1,2,\cdots,d)$ in Lemma \ref{xmlem1} (i). In Lemma \ref{xmlem1} (ii) and (iii), we shall also prove a recurrence formula for $\{(s;i_1,\cdots,i_d)\}$, which implies that $\Delta^L(k;\textbf{i})(\textbf{Y})=(m;1,2,\cdots,d)$ is expressed as a linear combination of
$\{(0;j_1,\cdots,j_d)| j_1,\cdots,j_d\in J,\ j_1<\cdots<j_d\}$. Note that if $(j_1,\cdots,j_d)=(m'+1,m'+2,\cdots,r,\ovl{d-r+m'},\ovl{d-r+m'-1},\cdots\ovl{1})$ (resp. $=(m'+1,m'+2,\cdots,m'+d)$), then $(0;j_1,\cdots,j_d)=1$ in the case $m'+d>r$ (resp. $m'+d\leq r$) by (\ref{ukeq1}), (\ref{ukeq2}) and (\ref{xmdef2}). If $(j_1,\cdots,j_d)$ is not as above, then $(0;j_1,\cdots,j_d)=0$. As a consequence of this calculation, we get Proposition \ref{pathlem}.

First, let us see the following lemma. We can verify it in the same way as (\ref{xmpro}).
\begin{lem}\label{xmpro2}For all $1\leq i,j\leq r-1$ and $Y\in \mathbb{C}^{\times}$, we have
\[
x_{-i}(Y)v_j=
\begin{cases}
Y^{-1}v_{i}+v_{i+1} & {\rm if}\ j=i, \\
Y v_{i+1} & {\rm if}\ j=i+1, \\
v_j & {\rm otherwise},
\end{cases}\q
x_{-i}(Y)v_{\ovl{j}}=
\begin{cases}
Y^{-1}v_{\ovl{i+1}}+v_{\ovl{i}} & {\rm if}\ j=i+1, \\
Y v_{\ovl{i}} & {\rm if}\ j=i, \\
v_{\ovl{j}} & {\rm otherwise}.
\end{cases}
\]
And for $1\leq i\leq r$, we also get
\[
x_{-i}(Y)v_r=
\begin{cases}
Y^{-2}v_{r}+Y^{-1}v_{0}+v_{\ovl{r}} & {\rm if}\ i=r, \\
Y v_{r-1} & {\rm if}\ i=r-1, \\
v_r & {\rm otherwise},
\end{cases}\q
x_{-i}(Y)v_{\ovl{r}}=
\begin{cases}
Y^{2}v_{\ovl{r}} & {\rm if}\ i=r, \\
Y^{-1}v_{\ovl{r}}+v_{\ovl{r-1}} & {\rm if}\ i=r-1, \\
v_{\ovl{r}} & {\rm otherwise},
\end{cases}
\]
\[
x_{-i}(Y)v_0=
\begin{cases}
v_{0}+2Y v_{\ovl{r}} & {\rm if}\ i=r, \\
v_0 & {\rm otherwise}.
\end{cases}
\]
\end{lem}

In the next lemma, we set $|l|=|\ovl{l}|=l$ for $1\leq l\leq r$.

\begin{lem}\label{xmlem1}
\begin{enumerate}
\item $\Delta^L(k;{\rm \bf{i}})(\textbf{Y})=(m;1,\cdots,d)$.

\item We suppose that $1\leq i_1<\cdots<i_{d}\leq \ovl{1}$ and $1\leq s\leq m$. Using the notation in Definition \ref{labeldef} $(i)$, we can write as follows:
\begin{eqnarray}\label{xmlem12}
& &(s;i_1,\cdots,i_d) \\
&=&\sum_{(j_1,\cdots,j_d)\in V}Q^{(m-s)}(i_1\rightarrow j_1)\cdots Q^{(m-s)}(i_d\rightarrow j_d)  \cdot(s-1;j_1,\cdots,j_d),  \nonumber
\end{eqnarray}
where $(j_1,\cdots,j_d)$ runs over $V:=\{(j_1,\cdots,j_d) |\ i_k\leq j_k\ (1\leq k\leq d),\ {\rm if\ }i_k\leq r-1\ {\rm then\ }j_k=i_k\ {\rm or}\ i_k+1,\ {\rm and\ }j_l\neq j_t\ (l\neq t) \}$. 

\item On the assumptions in $(ii)$, we can reduce the range $V$ of the sum in $(\ref{xmlem12})$ to $V':=\{(j_1,\cdots,j_d)\in V|\ {\rm If\ } \ovl{r}\leq j_l\ {\rm then}\ j_{l}\leq i_{l+1},\ {\rm and\ } j_{l}= i_{l+1}\ {\rm means}\ j_{l}= i_{l+1}=0 \}:$
\begin{eqnarray}\label{xmlem120}
& &(s;i_1,\cdots,i_d) \\
&=&\sum_{(j_1,\cdots,j_d)\in V'}a_{j_1,\cdots,j_d}Q^{(m-s)}(i_1\rightarrow j_1)\cdots Q^{(m-s)}(i_d\rightarrow j_d)  \cdot(s-1;j_1,\cdots,j_d),  \nonumber
\end{eqnarray}
where 
\[
a_{j_1,\cdots,j_d}=
\begin{cases}
\frac{1}{2} & {\rm if\ there\ exists\ some\ } l\ {\rm such\ that}\ j_l=i_{l+1}=0, \\
1 & {\rm otherwise}. 
\end{cases}
\]

\end{enumerate}

\end{lem}
\nd
{\sl Proof.}

(i) By Lemma \ref{xmpro2}, if $i>j$ $(i,j\in\{1,\cdots,r\})$, then we have $x_{-i}(Y)v_j=v_j$. Thus, we get
\begin{eqnarray*}
& &(m;1,\cdots,d):=\lan x^{(1)}_{-[1,r]} \cdots x^{(m-1)}_{-[1,r]}x^{(m)}_{-[1,r]}(v_1\wedge\cdots \wedge v_d),\ u_{\leq k}(v_1\wedge\cdots\wedge v_d)\ran \\
&=&\lan x^{(1)}_{-[1,r]} \cdots x^{(m-1)}_{-[1,r]}x_{-1}(Y_{m,1})\cdots x_{-d}(Y_{m,d})(v_1\wedge\cdots \wedge v_d),\ u_{\leq k}(v_1\wedge\cdots\wedge v_d)\ran\\
 &=&\lan x^L_{\textbf{i}}(\textbf{Y})(v_1\wedge\cdots \wedge v_d),\q u_{\leq k}(v_1\wedge\cdots\wedge v_d)\ran=\Delta^L(k;\textbf{i})(\textbf{Y}). 
\end{eqnarray*}

(ii) By Lemma \ref{xmpro2}, for $s$ $(1\leq s\leq m)$, we get
\begin{equation}\label{xmpro3}
x^{(s)}_{-[1,r]}v_i=
\begin{cases}
\frac{Y_{s,i-1}}{Y_{s,i}}v_i+v_{i+1} & {\rm if}\ 1\leq i\leq r-1, \\
\frac{Y_{s,r-1}}{Y^2_{s,r}}v_r+\frac{1}{Y_{s,r}}v_0+\sum^r_{j=1}\frac{1}{Y_{s,j-1}}v_{\ovl{j}}& {\rm if}\ i=r,
\end{cases}
\end{equation}
and
\begin{equation}\label{xmpro0}
x^{(s)}_{-[1,r]}v_{0}=v_0+2\sum^{r}_{j=1}\frac{Y_{s,r}}{Y_{s,j-1}}v_{\ovl{j}},
\end{equation}
\[
x^{(s)}_{-[1,r]}v_{\ovl{r}}=\frac{Y^2_{s,r}}{Y_{s,r-1}}v_{\ovl{r}}+\sum^{r-1}_{j=1}\frac{Y^2_{s,r}}{Y_{s,j-1}}v_{\ovl{j}},\ \ 
x^{(s)}_{-[1,r]}v_{\ovl{i}}=\sum^i_{j=1}\frac{Y_{s,i}}{Y_{s,j-1}}v_{\ovl{j}}\ \ (1\leq i\leq r-1),
\]
where we set $Y_{s,0}=1$. Note that the coefficient of $v_j$ in $x^{(s)}_{-[1,r]}v_i$ is equal to $Q^{(m-s)}(i\rightarrow j)$ $(1\leq i,j\leq \ovl{1})$ in Definition \ref{labeldef}. Thus,
\begin{eqnarray}\label{xmpro6} 
& &x^{(1)}_{-[1,r]} x^{(2)}_{-[1,r]}\cdots x^{(s-1)}_{-[1,r]}x^{(s)}_{-[1,r]}
(v_{i_1}\wedge\cdots\cdots \wedge v_{i_d}) \\
&=&x^{(1)}_{-[1,r]}\cdots x^{(s-1)}_{-[1,r]}
\sum_{j_1,\cdots,j_d} Q^{(m-s)}(i_1\rightarrow j_1)\cdots Q^{(m-s)}(i_d\rightarrow j_d)
(v_{j_1}\wedge\cdots\cdots \wedge v_{j_d}), \nonumber
\end{eqnarray}
where $(j_1,\cdots,j_d)$ runs over $V=\{(j_1,\cdots,j_d) |\ i_k\leq j_k\ (1\leq k\leq d),\ {\rm if\ }i_k\leq r-1\ {\rm then\ }j_k=i_k\ {\rm or}\ i_k+1,\ {\rm and\ }j_s\neq j_t\ (s\neq t)\}$. By pairing both sides in (\ref{xmpro6}) with $u_{\leq k}(v_1\wedge\cdots\wedge v_d)$, we obtain (\ref{xmlem12}).

(iii) We set $\hat{V}:=\{(j_1,\cdots,j_d)\in V|\ {\rm There\ exists\ a\ number\ }l\in[1,d-1]\ {\rm such\ that\ }\ovl{1}\geq j_l\geq i_{l+1}\geq0\}$. We define the map $\tau:\hat{V}\rightarrow \hat{V}$ as follows: We take $(j_1,\cdots,j_{\delta},j_{\delta+1},\cdots,j_d)\in \hat{V}$ and suppose that $j_{\delta}\in [1,r]$, $j_{\delta+1}\in \{\ovl{r},\cdots,\ovl{1}\}\cup\{0\}$. Let $l$ $(\delta+1\leq l\leq d-1)$ be the index such that $j_{\delta+1}<i_{\delta+2}, \cdots, j_{l-1}<i_{l}$ and $j_{l}\geq i_{l+1}$. Since $j_{l+1}\geq i_{l+1}$ by the definition of $V$, we have $(j_1,\cdots,j_{l+1},j_{l},\cdots,j_d)\in \hat{V}$. So, we define $\tau(j_1,\cdots,j_{l},j_{l+1},\cdots,j_d):= (j_1,\cdots,j_{l+1},j_{l},\cdots,j_d)$. We can easily see that $\tau^2=id_{\hat{V}}$. And then, we can verify that if $j_{l}<j_{l+1}$,

\vspace{3mm}

$Q(i_l\rightarrow j_l)Q(i_{l+1}\rightarrow j_{l+1})$
\[=
\begin{cases}
2Q(i_l\rightarrow j_{l+1})Q(i_{l+1}\rightarrow j_{l}) & {\rm if\ }i_l=r,\ i_{l+1}=0,\ {\rm and}\ j_{l}=0, \\
Q(i_l\rightarrow j_{l+1})Q(i_{l+1}\rightarrow j_{l}) & {\rm otherwise,}
\end{cases}
\]
by Definition \ref{labeldef} (i). Here, we have abbreviated $Q^{(m-s)}(i\rightarrow j)$ to $Q(i\rightarrow j)$. Furthermore, by (\ref{xmdef2}), it is easy to see
\[(s-1,j_1,\cdots,j_l,j_{l+1},\cdots,j_d)=-(s-1,j_1,\cdots,j_{l+1},j_l,\cdots,j_d). \]
Hence, if $j_l<j_{l+1}$ then the sum of $Q(i_l\rightarrow j_l)Q(i_{l+1}\rightarrow j_{l+1})(s-1;j_1,\cdots,j_l,j_{l+1},\cdots,j_d)$ and $Q(i_l\rightarrow j_{l+1})Q(i_{l+1}\rightarrow j_{l})(s-1;j_1,\cdots,j_{l+1},j_{l},\cdots,j_d)$ is equal to
\[
\begin{cases}
\frac{1}{2}Q(i_l\rightarrow j_{l})Q(i_{l+1}\rightarrow j_{l+1})(s-1;j_1,\cdots,j_l,j_{l+1},\cdots,j_d) & {\rm if\ }i_l=r,\ i_{l+1}=j_l=0, \\
0 & {\rm otherwise},
\end{cases}
\]
which implies that the partial sum of (\ref{xmlem120}) is reduced as follows: 
\begin{equation}\label{lem12pr1}
\sum_{(j_1,\cdots,j_d)\in\hat{V}}=
\begin{cases}
\frac{1}{2}\sum_{(j_1,\cdots,j_d)\in\hat{V},\ j_l=0} & {\rm if}\ i_l=r,\ i_{l+1}=0\ {\rm for}\ l\in[1,d-1],\\
0 & {\rm otherwise}.
\end{cases}
\end{equation}
We suppose that $i_l=r$, $i_{l+1}=0$ for some $l\in[1,d-1]$. In the set $\hat{V}_0:=\{(j_1,\cdots,j_d)\in\hat{V}| j_l=0,\ {\rm and\ there\ exists\ }\xi\ (\xi>l)\ {\rm such\ that\ }\ j_{\xi}\geq i_{\xi+1}\}$, we can also constitute the involution $\tau_0$ as follows: Let $(j_1,\cdots,j_d)\in\hat{V}_0$ and $j_{l+1}<i_{l+2}, \cdots, j_{\xi-1}<i_{\xi}$ and $j_{\xi}\geq i_{\xi+1}$. Then we define $\tau_0(j_1,\cdots,j_{\xi},j_{\xi+1},\cdots,j_d):= (j_1,\cdots,j_{\xi+1},j_{\xi},\cdots,j_d)$. By the same argument as above, we get 
\begin{equation}\label{lem12pr2}
\sum_{\hat{V}_0}=0.
\end{equation}
The reductions of range of sum (\ref{lem12pr1}) and (\ref{lem12pr2}) imply the equation (\ref{xmlem120}).
\qed

\vspace{3mm}

\nd
{\sl Proof of Proposition \ref{pathlem}.}

We suppose that $1\leq i_1<\cdots<i_d\leq\ovl{1}$ and $1\leq s\leq m$. By the definition of $V'$ in Lemma \ref{xmlem1}, we see that $(j_1,\cdots,j_d)\in V'$ if and only if the vertices vt$(s-1;j_1,\cdots,j_d)$ and vt$(s;i_1,\cdots,i_d)$ are connected (Definition \ref{connected-def}). Further, the coefficient of $(s-1;j_1,\cdots,j_d)$ in $(\ref{xmlem120})$ coincides with the label of the edge $e^{i_1,\cdots,i_d}_{j_1,\cdots,j_d}$ between vt$(s;i_1,\cdots,i_d)$ and vt$(s-1;j_1,\cdots,j_d)$ (Definition \ref{labeldef} (i)). Hence, we get
\begin{equation}\label{plpr1}
 (s;i_1,\cdots,i_d)=\sum_{(j_1,\cdots,j_d)} Q^{(m-s)}(e^{i_1,\cdots,i_d}_{j_1,\cdots,j_d}) \cdot(s-1;j_1,\cdots,j_d),
\end{equation}
where $(j_1,\cdots,j_d)$ runs over the set $\{(j_1,\cdots,j_d)|\ {\rm vt}(s-1;j_1,\cdots,j_d)$ and ${\rm vt}(s;i_1$,
$\cdots,i_d)\ {\rm are\ connected}\}$. Note that the conditions 
\[ {\rm if}\ j_l\in\{\ovl{r},\cdots,\ovl{1}\}\cup\{0\}\ {\rm then\ } j_l\leq i_{l+1},\]
\[ {\rm if}\ i_l\leq r-1\ {\rm then\ } j_l=i_l\ {\rm or}\ i_l+1, \]
and $i_{l+1}\leq j_{l+1}$ in $V'$ imply $j_l<j_{l+1}$, and we get $j_1<j_2<\cdots<j_d$. Using Lemma \ref{xmlem1} (iii), we obtain the followings in the same way as (\ref{plpr1}):
\begin{equation}\label{plpr2}
 (s-1;j_1,\cdots,j_d)=\sum_{(k_1,\cdots,k_d)} Q^{(m-s+1)}(e^{j_1,\cdots,j_d}_{k_1,\cdots,k_d}) \cdot(s-2;k_1,\cdots,k_d),
\end{equation}
where $(k_1,\cdots,k_d)$ runs over the set $\{(k_1,\cdots,k_d)|\ {\rm vt}(s-2;k_1,\cdots,k_d)$ and ${\rm vt}(s-1;j_1,\cdots,j_d)\ {\rm are\ connected}\}$ and $e^{j_1,\cdots,j_d}_{k_1,\cdots,k_d}$ is the edge between vertices ${\rm vt}(s-1;j_1,\cdots,j_d)$ and ${\rm vt}(s-2;k_1,\cdots,k_d)$. By (\ref{plpr1}), (\ref{plpr2}), $(s;i_1,\cdots,i_d)$ is a linear combination of $\{(s-2;k_1,\cdots,k_d)\}$, and the coefficient of $(s-2;k_1,\cdots,k_d)$ is as follows:
\[ \sum_{(j_1,\cdots,j_d)} Q^{(m-s)}(e^{i_1,\cdots,i_d}_{j_1,\cdots,j_d}) \cdot Q^{(m-s+1)}(e^{j_1,\cdots,j_d}_{k_1,\cdots,k_d}) \cdot(s-2;k_1,\cdots,k_d),\]
where $(j_1,\cdots,j_d)$ runs over the set $\{(j_1,\cdots,j_d)|\ {\rm vt}(s-1;j_1,\cdots,j_d)$ is connected to the vertices vt$(s;i_1,\cdots,i_d)$ and vt$(s-2;k_1,\cdots,k_d)\}$. The coefficient of $(s-2;k_1,\cdots,k_d)$ coincides with the label of subpath (Definition \ref{labeldef} (iii))
\[{\rm vt}(s;i_1,\cdots,i_d)\rightarrow {\rm vt}(s-1;j_1,\cdots,j_d)\rightarrow{\rm vt}(s-2;k_1,\cdots,k_d).\]

Repeating this argument, we see that $(s;i_1,\cdots,i_d)$ is a linear combination of $\{(0;l_1,\cdots,l_d)\}$ $(1\leq l_1<\cdots<l_d\leq \ovl{1})$. The coefficient of $(0;l_1,\cdots,l_d)$ is equal to the sum of labels of all subpaths from vt$(s;i_1,\cdots,i_d)$ to vt$(0;l_1,\cdots,l_d)$. In the case $m'+d>r$ (resp. $m'+d\leq r$), for $1\leq l_1<\cdots<l_d\leq \ovl{1}$, if $(l_1,\cdots,l_d)=(m'+1,m'+2,\cdots,r,\ovl{d-r+m'},\cdots,\ovl{2},\ovl{1})$ (resp. $=(m'+1,m'+2,\cdots,m'+d)$), then we obtain $(0;l_1,\cdots,l_d)=1$ by (\ref{ukeq1}), (\ref{ukeq2}) and (\ref{xmdef2}). If $(l_1,\cdots,l_d)$ is not as above, we obtain $(0;l_1,\cdots,l_d)=0$. Therefore, we see that $(s;i_1,\cdots,i_d)$ is equal to the sum of labels of subpaths from vt$(s;i_1,\cdots,i_d)$ to vt$(0;m'+1,m'+2,\cdots,r,\ovl{d-r+m'},\cdots,\ovl{2},\ovl{1})$ (resp. vt$(0;m'+1,m'+2,\cdots,m'+d))$.

In particular, $\Delta^L(k;\textbf{i})(\textbf{Y})=(m;1,2,\cdots,d)$ is equal to the sum of labels of paths in $X_d(m,m')$, which means $\Delta^L(k;\textbf{i})(\textbf{Y})=\sum_{p\in
X_d(m,m')} Q(p)$. 
\qed

\begin{ex}\label{pathex2}
Let us consider the same setting as in Example \ref{pathex3}, {\it i.e.},
$r=3$, $u=s_1s_2s_3s_1s_2s_3s_1s_2$, $v=e$, $k=5$, 
${\rm \bf{i}}=(1,2,3,1,2,3,1,2)$, 
$m=3$, $m'=2$ and $d=2$.
Therefore,  by Example \ref{pathex}, we obtain
\begin{eqnarray*}
\Delta^L(5;{\rm \bf{i}})(\textbf{Y})&=&\sum^{15}_{i=1} Q(p_i)\\
&=&
\frac{1}{Y_{3,2}}+\frac{Y_{2,2}}{Y_{3,1}Y^2_{2,3}}+2\frac{Y_{1,3}}{Y_{3,1}Y_{2,3}}+\frac{Y^2_{1,3}}{Y_{3,1}Y_{2,2}}
+\frac{Y_{1,2}}{Y_{3,1}Y_{2,1}}+\frac{Y_{1,1}}{Y_{3,1}}\\
\q &+&\frac{Y_{2,1}}{Y^2_{2,3}}+2\frac{Y_{2,1}Y_{1,3}}{Y_{2,2}Y_{2,3}}+\frac{Y_{2,1}Y^2_{1,3}}{Y^2_{2,2}}+2\frac{Y_{1,2}}{Y_{2,2}}+
\frac{Y_{2,1}Y_{1,1}}{Y_{2,2}}+2\frac{Y_{1,2}}{Y_{2,3}Y_{1,3}}\\
&+&\frac{Y^2_{1,2}}{Y_{2,1}Y^2_{1,3}}+\frac{Y_{1,1}Y_{1,2}}{Y^2_{1,3}}. 
\end{eqnarray*}
We find that this coincides with the explicit form of 
$\Delta^L(5;{\rm \bf{i}})(\textbf{Y})$ in Example \ref{pathex3}.
\end{ex}

\begin{rem}\label{coinArem}
We suppose that $m'+d\leq r$.
\begin{enumerate}
\item[$(1)$] Definition \ref{pathdef} shows that the set $X_d(m,m')$ is constituted by paths $p$
\begin{multline*}p={\rm vt}(m;a^{(0)}_1,\cdots,a^{(0)}_d)\rightarrow{\rm vt}(m-1;a^{(1)}_1,\cdots,a^{(1)}_d)\rightarrow\\
\cdots\rightarrow{\rm vt}(1;a^{(m-1)}_1,\cdots,a^{(m-1)}_d)\rightarrow{\rm vt}(0;a^{(m)}_1,\cdots,a^{(m)}_d)
\end{multline*}
which satisfy the following conditions: For $0\leq s\leq m$,
\begin{enumerate}
\item[$(i)$] $a^{(s)}_{\zeta}\in\{1,\cdots,r\}$ $(1\leq \zeta\leq d)$, 
\item[$(ii)$] $a^{(s)}_{1}<a^{(s)}_{2}<\cdots<a^{(s)}_{d}$, 
\item[$(iii)$] $a^{(s+1)}_{\zeta}=a^{(s)}_{\zeta}$ or $a^{(s)}_{\zeta}+1$,
\item[$(iv)$] $(a^{(0)}_1,a^{(0)}_2,\cdots,a^{(0)}_d)=(1,2,\cdots,d)$, 

$(a^{(m)}_1,\cdots,a^{(m)}_d)=(m'+1,m'+2,\cdots,m'+d)$.
\end{enumerate}
\item[$(2)$] By Definition \ref{labeldef}, the label $Q^{(s)}(p)$ of the edge ${\rm vt}(m-s;a^{(s)}_1,a^{(s)}_2,\cdots,a^{(s)}_d)\rightarrow{\rm vt}(m-s-1;a^{(s+1)}_1,a^{(s+1)}_2,\cdots,a^{(s+1)}_d)$ is as follows:
\[
Q^{(s)}(p):=\frac{Y_{m-s,a^{(s+1)}_1-1}}{Y_{m-s,a^{(s)}_1}}\cdots \frac{Y_{m-s,a^{(s+1)}_d-1}}{Y_{m-s,a^{(s)}_d}}.
\]
\item[$(3)$] For $G_A=SL_{r+1}(\mathbb{C})$, let $B_A$ and $(B_-)_A$ be two opposite Borel subgroups in $G_A$, $N_A\subset B_A$ and $(N_-)_A\subset (B_-)_A$ their unipotent radicals, and $W_A$ be the Weyl group of $G_A$. We define a reduced double Bruhat cell as $L^{u,v}_A:= (N_A\cdot \ovl{u}\cdot N_A) \cap ((B_-)_A\cdot \ovl{v}\cdot (B_-)_A)$.
We set $u,v\in W_A$ and their reduced word ${\rm \bf{i}}_A$ as
\[ u=\underbrace{s_1\cdots s_r}_{1\ {\rm st\ cycle}}\underbrace{s_1\cdots s_{r-1}}_{2\ {\rm nd\ cycle}}\cdots \underbrace{s_1\cdots s_{i_n}}_{m\ {\rm th\ cycle}}, \q v=e, \]
\[ {\rm \bf{i}}_A=(\underbrace{1,\cdots, r}_{1\ {\rm st\ cycle}},\underbrace{1,\cdots, (r-1)}_{2\ {\rm nd\ cycle}},\cdots,\underbrace{1,\cdots ,i_n}_{m\ {\rm th\ cycle}}), \]
where $n=l(u)$ and $1\leq i_n\leq r-m+1$. Let $i_k$ be the $k$-th index of  ${\rm \bf{i}}_A$ from the left, and belong to $m'$-th cycle. Using Theorem \ref{fp2}, we can define $\Delta^{L_A}(k;{\rm \bf{i}}_A)(\textbf{Y}_A):=(\Delta(k;{\rm \bf{i}}_A)\circ
 x^{L_A}_{{\rm \bf{i}}_A})(\textbf{Y}_A)$ in the same way as Definition \ref{gendef}, where 
\[
 \textbf{Y}_A:=(Y_{1,1},Y_{1,2},\cdots,Y_{1,r},Y_{2,1},Y_{2,2},\cdots,Y_{2,r-1},\cdots,Y_{m,1},\cdots,Y_{m,i_n})\in  (\mathbb{C}^{\times})^{n},
\]
and the map $x^{L_A}_{{\rm \bf{i}}_A}:(\mathbb{C}^{\times})^{n} \overset{\sim}{\hookrightarrow} L^{u,v}_A$ is defined as in Theorem \ref{fp2}.

Then, we already have seen in \cite{KaN} that 
$\Delta^{L_A}(k;{\rm \bf{i}}_A)(\textbf{Y}_A)=\sum_{p\in X_d(m,m')}Q(p)$, 
where $X_d(m,m')$ and the label $Q(p)=\prod^{m-1}_{s=0}Q^{(s)}(p)$ 
is the one we have seen in $(1)$ and $(2)$. 
Therefore, it follows from Proposition \ref{pathlem} that 
if $m'+d\leq r$, then $\Delta^L(k;{\rm \bf{i}})(\textbf{Y})$ 
coincides with $\Delta^{L_A}(k;{\rm \bf{i}}_A)(\textbf{Y}_A)$.
\end{enumerate}
\end{rem}

\subsection{The properties of paths in $X_d(m,m')$}

In this subsection, we shall see some lemmas on $X_d(m,m')$. By Remark \ref{coinArem}, we suppose that $m'+d>r$. We fix a path $p\in X_d(m,m')$
\begin{multline}\label{fixpath}
p={\rm vt}(m;a^{(0)}_1,\cdots,a^{(0)}_d)\rightarrow
\cdots\rightarrow{\rm vt}(2;a^{(m-2)}_1,\cdots,a^{(m-2)}_d)  \\
\rightarrow{\rm vt}(1;a^{(m-1)}_1,\cdots,a^{(m-1)}_d)\rightarrow{\rm vt}(0;a^{(m)}_1,\cdots,a^{(m)}_d).
\end{multline}

\begin{lem}\label{cancellem}
For $p\in X_d(m,m')$ in $(\ref{fixpath})$ and $i$ $(r-m'+1\leq i\leq d)$, we obtain 
\begin{equation}\label{fixedeq}
a^{(m)}_i=a^{(m-1)}_i=\cdots=a^{(m-i+r-m'+1)}_i=\ovl{d-i+1}.
\end{equation}
\end{lem}
\nd
{\sl Proof.}

By Definition \ref{pathdef} (iv), we get $a^{(m)}_{r-m'+1}=\ovl{d-r+m'}$, and by Definition \ref{pathdef} (v), we also get $\ovl{d-r+m'}=a^{(m)}_{r-m'+1}<a^{(m-1)}_{r-m'+2}\leq \ovl{1}$. Using Definition \ref{pathdef} (v) repeatedly, we obtain $\ovl{d-r+m'}=a^{(m)}_{r-m'+1}<a^{(m-1)}_{r-m'+2}<a^{(m-2)}_{r-m'+3}<\cdots<a^{(m-d+r-m'+1)}_{d}\leq \ovl{1}$, which means
\[a^{(m-i+r-m'+1)}_{i}=\ovl{d-i+1}\ \ (r-m'+1\leq i\leq d). \]
It follows from (\ref{iseqine}) and Definition \ref{pathdef} (iv) that $\ovl{d-i+1}=a^{(m-i+r-m'+1)}_{i}\leq a^{(m-i+r-m'+2)}_{i}\leq\cdots\leq a^{(m-1)}_{i}\leq a^{(m)}_{i}=\ovl{d-i+1}$, which yields (\ref{fixedeq}). \qed

\vspace{3mm}

By this lemma, we get $a^{(s)}_i=\ovl{d-i+1}$ for $s$ $(m-i+r-m'+1\leq s\leq m)$. In the next lemma, we see the properties for $a^{(s)}_i$ $(0\leq s\leq m-i+r-m')$.

\begin{lem}\label{tab-length}
For $i$ $(1\leq i\leq d)$ and $p\in X_d(m,m')$, let $a^{(0)}_{i}\rightarrow a^{(1)}_{i}\rightarrow a^{(2)}_{i}\rightarrow \cdots\rightarrow a^{(m)}_{i}$
be the $i$-sequence of the path $p$ $(Definition\ \ref{iseq})$.
\begin{enumerate}
\item In the case $i\leq r-m'$,
\[   \#\{0\leq s\leq m-1|\ 1\leq a^{(s)}_i\leq r,\ {\rm and}\ a^{(s)}_i=a^{(s+1)}_i\}=m-m'. \]
\item In the case $i> r-m'$,
\begin{multline*}
\#\{0\leq s\leq m-i+r-m'|\ 1\leq a^{(s)}_i\leq r,\ {\rm and}\ a^{(s)}_i=a^{(s+1)}_i\}+\\
\#\{0\leq s\leq m-i+r-m'|\ a^{(s)}_i\in\{\ovl{r},\cdots,\ovl{1}\}\cup\{0\}\}=m-m'. 
\end{multline*}
\end{enumerate}
\end{lem}
\nd
{\sl Proof.}

(i) In the case $i\leq r-m'$, Definition \ref{pathdef} (iii), (iv) and (\ref{iseqine}) show that 
\begin{equation}\label{tab-length1}
 i=a^{(0)}_i\leq a^{(1)}_i\leq\cdots\leq a^{(m)}_i=m'+i,\ \ a^{(s+1)}_i=a^{(s)}_i\ {\rm or}\ a^{(s)}_i+1. 
\end{equation}
In particular, we get $1\leq a^{(s)}_i\leq r$ for $1\leq s\leq m$.
By (\ref{tab-length1}), we obtain
\[ \#\{0\leq s\leq m-1|\ a^{(s+1)}_i=a^{(s)}_i+1\}=m',
\]
which implies $\#\{0\leq s\leq m-1|\ a^{(s)}_i=a^{(s+1)}_i\}=m-m'$. 

(ii) In the case $i>r-m'$, by (\ref{iseqine}), we have
\[
 i=a^{(0)}_i\leq a^{(1)}_i\leq\cdots\leq a^{(m-i+r-m')}_i\leq \ovl{1}. 
\]
We suppose that 
\begin{equation}\label{tab-length2}
i= a^{(0)}_i\leq a^{(1)}_i\leq \cdots\leq a^{(l)}_i\leq r,\ {\rm and}\ \ 0\leq a^{(l+1)}_i\leq \cdots\leq a^{(m-i+r-m')}_i\leq\ovl{1},
\end{equation}
for some $l$ $(1\leq l\leq m-i+r-m')$. Definition \ref{pathdef} (iii) implies that $a^{(s+1)}_i=a^{(s)}_i$ or $a^{(s)}_i+1$ $(1\leq s\leq l-1)$ and $a^{(l)}_i=r$. Therefore, 
\[ i=a^{(0)}_i\leq a^{(1)}_i\leq\cdots\leq a^{(l)}_i=r,\qq a^{(s+1)}_i=a^{(s)}_i\ {\rm or}\ a^{(s)}_i+1. \]
So we have $\#\{0\leq s\leq l-1 |\ a^{(s+1)}_i=a^{(s)}_i \}=l-(r-i)$ in the same way as (i).

On the other hand, the assumption $0\leq a^{(l+1)}_i\leq\cdots\leq a^{(m-i+r-m')}_i\leq \ovl{1}$ in (\ref{tab-length2}) means that $\{l+1\leq s\leq m-i+r-m' |\ a^{(s)}_i\in \{\ovl{r},\cdots,\ovl{1}\}\cup\{0\}\}=\{l+1,l+2,\cdots, m-i+r-m'\}$. Hence, $\#\{1\leq s\leq l-1 |\ a^{(s+1)}_i=a^{(s)}_i \}+\#\{l+1\leq s\leq m-i+r-m' |\ a^{(s)}_i\in \{\ovl{r},\cdots,\ovl{1}\}\cup\{0\}\}=(l-(r-i))+(m-i+r-m'-l)=m-m'$. \qed

\vspace{3mm}

By this lemma, we define $l^{(s)}_i\in\{0,1,\cdots,m\}$ $(1\leq i\leq d,\ 1\leq s\leq m-m')$ for the path $p\in X_d(m,m')$ in $(\ref{fixpath})$ as follows:  For $i\leq r-m'$, we set $\{l^{(s)}_i\}_{1\leq s\leq m-m'}$ $(l^{(1)}_i<\cdots<l^{(m-m')}_i)$ as
\begin{equation}\label{ldef}
\{l^{(1)}_i,\ l^{(2)}_i,\cdots,l^{(m-m')}_i\}:=\{s|a^{(s)}_i=a^{(s+1)}_i,\ \ 0\leq s\leq m-1 \}.
\end{equation}
For $i>r-m'$, we set $\{l^{(s)}_i\}_{1\leq s\leq m-m'}$ $(l^{(1)}_i<\cdots<l^{(m-m')}_i)$ as
\begin{eqnarray}\label{ldef2}
& &\{l^{(1)}_i,\ l^{(2)}_i,\cdots,l^{(m-m')}_i\}\nonumber \\
&:=&\{s|\ 1\leq a^{(s)}_i\leq r,\ a^{(s)}_i=a^{(s+1)}_i,\ 0\leq s\leq m-i+r-m' \}\nonumber \\
&\cup&\ \{s|\ a^{(s)}_i\in \{\ovl{r},\cdots,\ovl{1}\}\cup\{0\}, \ 0\leq s\leq m-i+r-m' \}. \qq \q
\end{eqnarray}
We also set $k^{(s)}_i\in \{j,\ovl{j}|\ 1\leq j\leq r\}\cup\{0\}$ $(1\leq i\leq d,\ 1\leq s\leq m-m')$ as
\begin{equation}\label{kdef}
k^{(s)}_i:=a^{(l^{(s)}_i)}_i.
\end{equation}
Using (\ref{iseqine}) and $l^{(1)}_i<\cdots<l^{(m-m')}_i$, we obtain
\begin{equation}\label{kpro1}
\begin{cases}
i\leq k^{(1)}_i\leq \cdots\leq k^{(m-m')}_i\leq m'+i & {\rm if\ }1\leq i\leq r-m', \\
i\leq k^{(1)}_i\leq \cdots\leq k^{(m-m')}_i\leq \ovl{d-i+1} & {\rm if\ }r-m'< i\leq d. \\
\end{cases}
\end{equation}
For $i$ $(1\leq i\leq d)$, by (\ref{kpro1}), we can define the number $\delta_i$ $(0\leq\delta_i\leq m-m',\ \delta_{d}\leq\cdots\leq\delta_2\leq\delta_1)$ as
\begin{equation}\label{deldef}
1\leq k^{(1)}_i\leq\cdots\leq k^{(\delta_i)}_i\leq r,\q 0\leq k^{(\delta_i+1)}_i\leq\cdots\leq k^{(m-m')}_i\leq\ovl{1}.
\end{equation}

\begin{lem}\label{kpro2}
\begin{enumerate}
\item For $1\leq i\leq d$,
\[ l^{(s)}_{i}=
\begin{cases}
k^{(s)}_{i}+s-i-1 & {\rm if}\  k^{(s)}_{i}\in\{j|1\leq j\leq r\}, \\
s-i+r & {\rm if}\ k^{(s)}_{i}\in\{\ovl{j}|1\leq j\leq r\}\cup\{0\}.
\end{cases}
\]
\item For $1\leq s\leq m-m'$ and $1\leq i\leq d-1$, if $k^{(s)}_{i}\in\{j|1\leq j\leq r\}$, then we have $k^{(s)}_{i}<k^{(s)}_{i+1}$, and $l^{(s)}_i\leq l^{(s)}_{i+1}$.

For $1\leq i\leq d-1$, if $k^{(s)}_{i}\in\{\ovl{j}|1\leq j\leq r\}\cup\{0\}$, then we have
$k^{(s)}_{i}\leq k^{(s)}_{i+1}$, $l^{(s)}_{i}=l^{(s)}_{i+1}+1$,
and $k^{(s)}_{i}=k^{(s)}_{i+1}$ if and only if $k^{(s)}_{i}=k^{(s)}_{i+1}=0$. 
\item For $1\leq s\leq m-m'-1$ and $1\leq i\leq d-1$, we have $k^{(s)}_{i}<k^{(s+1)}_{i+1}$.
\end{enumerate}
\end{lem}
\nd
{\sl Proof.}

(i) We suppose that $k^{(s)}_{i}\in \{j|1\leq j\leq r\}$. The definition of $l^{(s)}_i$ in (\ref{ldef}) means that the path $p$ has the following $i$-sequence (Definition \ref{iseq}):
\begin{flushleft}
$a^{(0)}_i=i,\ a^{(1)}_i=i+1,\ a^{(2)}_i=i+2,\cdots,a^{(l^{(1)}_i)}_i=i+l^{(1)}_i$,
 
$ a^{(l^{(1)}_i+1)}_i=i+l^{(1)}_i,\ a^{(l^{(1)}_i+2)}_i=i+l^{(1)}_i+1,\cdots, a^{(l^{(2)}_i)}_i=i+l^{(2)}_i-1, $

$ a^{(l^{(2)}_i+1)}_i=i+l^{(2)}_i-1,\ a^{(l^{(2)}_i+2)}_i=i+l^{(2)}_i,\cdots, a^{(l^{(3)}_i)}_i=i+l^{(3)}_i-2, $
\end{flushleft}
\begin{equation}\label{jlist} \vdots \end{equation}
\[ a^{(l^{(s-1)}_i+1)}_i=i+l^{(s-1)}_i-s+2,\ a^{(l^{(s-1)}_i+2)}_i=i+l^{(s-1)}_i-s+3,\cdots, a^{(l^{(s)}_i)}_i=i+l^{(s)}_i-s+1,\]
$a^{(l^{(s)}_i+1)}_i=i+l^{(s)}_i-s+1,\ a^{(l^{(s)}_i+2)}_i=i+l^{(s)}_i-s+2,\cdots$.

\vspace{2mm}

Hence we have $k^{(s)}_{i}=a^{(l^{(s)}_{i})}_i=i+l^{(s)}_{i}-s+1$, which implies $l^{(s)}_{i}=k^{(s)}_{i}+s-i-1$.

Next, we suppose that $a^{(l^{(s)}_{i})}_i=k^{(s)}_{i}\in \{\ovl{j}|1\leq j\leq r\}\cup\{0\}$. Using (\ref{iseqine}), we get $a^{(l^{(s)}_{i})}_i\leq a^{(l^{(s)}_{i}+1)}_i\leq\cdots\leq a^{(m-i+r-m')}_i$, so $a^{(\zeta)}_i\in\{\ovl{j}|1\leq j\leq r\}\cup\{0\}$ $(l^{(s)}_{i}\leq \zeta \leq m-i+r-m')$. Thus, by the definition (\ref{ldef2}) of $l^{(s)}_i$, we obtain $l^{(m-m')}_i=m-i+r-m'$, $l^{(m-m'-1)}_i=m-i+r-m'-1,l^{(m-m'-2)}_i=m-i+r-m'-2,\cdots, l^{(\xi)}_i=\xi-i+r$ $(s\leq\xi\leq m-m')$. In particular, we get
\begin{equation}\label{kpro2pr10}
 l^{(s)}_i=s-i+r.
\end{equation}

(ii) We suppose that $k^{(s)}_{i}\in \{j|1\leq j\leq r\}$. If $k^{(s)}_{i+1}\in \{\ovl{j}|1\leq j\leq r\}\cup\{0\}$, then clearly $k^{(s)}_{i}<k^{(s)}_{i+1}$ and it follows from (i) that $l^{(s)}_{i}\leq l^{(s)}_{i+1}$. So we may assume that $k^{(s)}_{i+1}\in \{j|1\leq j\leq r\}$. By Definition \ref{pathdef} (ii) and the definition (\ref{ldef2}) of $l^{(s)}_{i+1}$, we have $a^{(l^{(s)}_{i+1}+1)}_i<a^{(l^{(s)}_{i+1}+1)}_{i+1}=a^{(l^{(s)}_{i+1})}_{i+1}=k^{(s)}_{i+1}\leq r$. Therefore, the inequality (\ref{iseqine}) implies
\begin{equation}\label{kpro2pr3ano}
i=a^{(0)}_{i}\leq a^{(1)}_{i}\leq\cdots\leq a^{(l^{(s)}_{i+1})}_{i}\leq a^{(l^{(s)}_{i+1}+1)}_{i}<r, 
\end{equation}
\[ a^{(\zeta)}_{i}=a^{(\zeta-1)}_{i}\ {\rm or}\ a^{(\zeta-1)}_{i}+1\ (1\leq\zeta\leq l^{(s)}_{i+1}+1). \]
We obtain
\begin{equation}\label{kpro2pr4ano}
l^{(s)}_{i+1}+1-s\geq\#\{\zeta| a^{(\zeta)}_{i}=a^{(\zeta-1)}_{i}+1,\ 1\leq\zeta\leq l^{(s)}_{i+1}+1\},
\end{equation}
otherwise, it follows from (\ref{kpro2pr3ano}) and (i) that $a^{(l^{(s)}_{i+1}+1)}_i> i+l^{(s)}_{i+1}+1-s=k^{(s)}_{i+1}-1=a^{(l^{(s)}_{i+1})}_{i+1}-1$, and hence $a^{(l^{(s)}_{i+1}+1)}_i\geq a^{(l^{(s)}_{i+1})}_{i+1}=a^{(l^{(s)}_{i+1}+1)}_{i+1}$, which contradicts Definition \ref{pathdef} (ii).

The inequality (\ref{kpro2pr4ano}) yields that
\begin{equation}\label{kpro2pr5ano}
s\leq\#\{\zeta| a^{(\zeta)}_{i}=a^{(\zeta-1)}_{i},\ 1\leq\zeta\leq l^{(s)}_{i+1}+1\}.
\end{equation}
On the other hand, the definition of $l^{(s)}_i$ implies $a^{(l^{(s)}_i+1)}_i=a^{(l^{(s)}_i)}_i=k^{(s)}_{i}\in\{j|1\leq j\leq r\}$. The inequality (\ref{iseqine}) shows 
\[
i=a^{(0)}_{i}\leq a^{(1)}_{i}\leq \cdots\leq a^{(l^{(s)}_{i})}_{i}=a^{(l^{(s)}_{i}+1)}_{i}=k^{(s)}_{i},
\] 
\[ a^{(\zeta)}_{i}=a^{(\zeta-1)}_{i}\ {\rm or}\ a^{(\zeta-1)}_{i}+1\ (1\leq\zeta\leq l^{(s)}_{i}+1),\]
and
\begin{equation}\label{kpro2pr6ano}
s=\#\{\zeta| a^{(\zeta)}_{i}=a^{(\zeta-1)}_{i},\ 1\leq\zeta\leq l^{(s)}_{i}+1\}.
\end{equation}
Since $a^{(l^{(s)}_i)}_i=a^{(l^{(s)}_i+1)}_i$, the equation (\ref{kpro2pr6ano}) means 
\begin{equation}\label{kpro2pr6ano2}
s-1=\#\{\zeta| a^{(\zeta)}_{i}=a^{(\zeta-1)}_{i},\ 1\leq\zeta\leq l^{(s)}_{i}\}.
\end{equation}
Thus, by (\ref{kpro2pr5ano}) and (\ref{kpro2pr6ano2}), we have $l^{(s)}_{i}< l^{(s)}_{i+1}+1$, and hence $l^{(s)}_{i}\leq l^{(s)}_{i+1}$, which yields $k^{(s)}_i<k^{(s)}_{i+1}$ since $k^{(s)}_i=i+l^{(s)}_i-s+1<i+l^{(s)}_{i+1}-s+2=(i+1)+l^{(s)}_{i+1}-s+1=k^{(s)}_{i+1}$.

Next, we suppose that $a^{(l^{(s)}_i)}_{i}=k^{(s)}_{i}\in \{\ovl{j}|1\leq j\leq r\}\cup\{0\}$. It follows from Definition \ref{pathdef} (v) that $a^{(l^{(s)}_i-1)}_{i+1}\in \{\ovl{j}|1\leq j\leq r\}\cup\{0\}$. Since $a^{(l^{(s)}_i-1)}_{i+1}\leq a^{(l^{(s)}_i)}_{i+1}\leq \cdots\leq a^{(m-(i+1)+r-m')}_{i+1}$, we get $a^{(\zeta)}_{i+1}\in\{\ovl{j}|1\leq j\leq r\}\cup\{0\}$ $(l^{(s)}_i-1\leq\zeta\leq m-(i+1)+r-m')$. By the same way as in (\ref{kpro2pr10}), we can verify $l^{(s)}_{i+1}=s-(i+1)+r$. Therefore, it follows from (\ref{kpro2pr10}) that $l^{(s)}_{i}=l^{(s)}_{i+1}+1$. Further, $k^{(s)}_i=a^{(l^{(s)}_i)}_i\leq a^{(l^{(s)}_i-1)}_{i+1}=a^{(l^{(s)}_{i+1})}_{i+1}=k^{(s)}_{i+1}$, and $k^{(s)}_i=k^{(s)}_{i+1}$ if and only if $k^{(s)}_i=k^{(s)}_{i+1}=0$ by Definition \ref{pathdef} (v). 

(iii) It is clear that if $k^{(s)}_{i}\neq0$ then $k^{(s)}_{i}<k^{(s+1)}_{i+1}$ by (ii) and (\ref{kpro1}). If $k^{(s)}_{i}=0$ then we get $0=k^{(s)}_{i}\leq k^{(s+1)}_{i+1}$. Assuming that $0=k^{(s)}_{i}=k^{(s+1)}_{i+1}$, we obtain $a^{(s-i+r)}_{i}=a^{(s-i+r)}_{i+1}=0$ by (i). It contradicts Definition \ref{pathdef} (ii). Thus, we get $k^{(s)}_{i}<k^{(s+1)}_{i+1}$.
\qed

\subsection{The proof of Theorem \ref{thm1}}

First, we see the following lemma. Let us recall the definition (\ref{bbbar}) of $B(l,k)$.
\begin{lem}\label{thm1lem}
For $p\in X_d(m,m')$ in $(\ref{fixpath})$, we set $l^{(s)}_i$, $k^{(s)}_i$ and $\delta_i$ as in $(\ref{ldef})$, $(\ref{ldef2})$, $(\ref{kdef})$ and $(\ref{deldef})$. Then we have
\begin{multline}\label{thm1lemclaim}
Q(p)= 2^{C[\{k^{(s)}_i\}]}\prod^{d}_{i=1} B(m-l^{(1)}_i,k^{(1)}_i)\cdots B(m-l^{(\delta_i)}_i,k^{(\delta_i)}_i)\\
\cdot B(m-l^{(\delta_i+1)}_i,k^{(\delta_i+1)}_i)\cdots B(m-l^{(m-m')}_i,k^{(m-m')}_i),
\end{multline}
where 
\begin{multline}\label{Cpdef}
C[\{k^{(s)}_i\}]=C[\{k^{(s)}_i\}_{1\leq i\leq d,\ 1\leq s\leq m-m'}]:=\\
\#\{j\in[1,d-1]|\ k^{(t)}_{j+1}=0,\ k^{(t+1)}_{j+1}\neq 0\ {\rm and}\ k^{(t+1)}_j\neq0\ {\rm for\ some}\ t\in[1,m-m']\}.
\end{multline}
\end{lem}
\nd
{\sl Proof.}

Recall that the numbers $\delta_i$ $(1\leq i\leq d)$ are defined as $(\ref{deldef})$, which implies that the $i$-sequence (Definition \ref{iseq}) of the path $p$ satisfies the following inequality in the order $(\ref{B-order})$:
\begin{equation}\label{thm1lem-rel1}
1\leq a^{(0)}_i\leq a^{(1)}_i\leq\cdots \leq a^{(l^{(\delta_i+1)}_i-1)}_i=r<0\leq a^{(l^{(\delta_i+1)}_i)}_i\leq \cdots\leq a^{(m)}_i\leq\ovl{1}.
\end{equation}

Next, let us also recall the definition of the labels of the path $p$ in Definition \ref{labeldef}:
\[ Q^{(s)}(p)=c^{(s)}(p)\prod^{d}_{i=1} Q^{(s)}(a^{(s)}_i\rightarrow a^{(s+1)}_i), \]
where 
\[
c^{(s)}(p)=
\begin{cases}
\frac{1}{2} & {\rm if\ there\ exists\ some\ } i\ {\rm such\ that}\ a^{(s)}_{i+1}=a^{(s+1)}_{i}=0, \\
1 & {\rm otherwise},
\end{cases}
\]
and
\[ Q(p)=\prod^{m-1}_{s=0}c^{(s)}(p)\prod^{d}_{i=1} Q^{(s)}(a^{(s)}_i\rightarrow a^{(s+1)}_i).\]
To calculate $Q(p)$, let us divide the range of product $\prod^{m-1}_{s=0}$ as follows:
\[
\prod^{l^{(\delta_i+1)}_{i}-2}_{s=0},\qq 
\prod^{l^{(m-m')}_i}_{s=l^{(\delta_i+1)}_i-1}\qq {\rm and}\qq
\prod^{m-1}_{s=l^{(m-m')}_i+1},
\]
where in the case $\delta_i=m-m'$, we set
\begin{equation}\label{dmdmdm}
l^{(m-m'+1)}_i:=l^{(m-m')}_i+2.
\end{equation}
 
First, let us consider the first range $0\leq s\leq l^{(\delta_i+1)}_i-2$ of the product. In this range, the inequality (\ref{thm1lem-rel1}) means $a^{(s)}_i\leq a^{(s+1)}_i\leq r$. By Definition \ref{labeldef} (i),
\[ Q^{(s)}(a^{(s)}_i\rightarrow a^{(s+1)}_i)=
\begin{cases}
\frac{Y_{m-s,a^{(s)}_i-1}}{Y_{m-s,a^{(s)}_i}}=B(m-s,a^{(s)}_i)  & {\rm if}\ a^{(s)}_i=a^{(s+1)}_i\leq r-1, \\
\frac{Y_{m-s,r-1}}{Y^2_{m-s,r}}=B(m-s,r)  & {\rm if}\ a^{(s)}_i=a^{(s+1)}_i=r, \\
1 & {\rm if}\ a^{(s+1)}_i=a^{(s)}_i+1,
\end{cases}
 \]
which implies that
\begin{equation*}
\prod^{l^{(\delta_i+1)}_{i}-2}_{s=0}\left(Q^{(s)}(a^{(s)}_i\rightarrow a^{(s+1)}_i)\right)
=\prod^{\delta_i}_{\zeta=1}B(m-l^{(\zeta)}_i,k^{(\zeta)}_i),
\end{equation*}
where we used the notation in (\ref{ldef2}) and $k^{(\zeta)}_i=a^{(l^{(\zeta)}_i)}_i$. Hence,
\begin{equation}\label{thm1lempr2}
\prod^{l^{(\delta_i+1)}_{i}-2}_{s=0}c^{(s)}(p)\prod^{d}_{i=1} \left(Q^{(s)}(a^{(s)}_i\rightarrow a^{(s+1)}_i)\right)
=\prod^{d}_{i=1} \prod^{\delta_i}_{\zeta=1}B(m-l^{(\zeta)}_i,k^{(\zeta)}_i).
\end{equation}

Next, we consider the second range $l^{(\delta_i+1)}_i-1\leq s\leq l^{(m-m')}_i $ of the product. In this range, the inequality (\ref{thm1lem-rel1}) means $a^{(s+1)}_i\in\{\ovl{j}|1\leq j\leq r\}\cup\{0\}$. If $r-m'\geq i$ then $r\geq m'+i=a^{(m)}_i\geq\cdots\geq a^{(1)}_i\geq a^{(0)}_i$, which implies $\delta_i=m-m'$, and $\prod^{l^{(m-m')}_i}_{s=l^{(\delta_i+1)}_i-1}\left(Q^{(s)}(a^{(s)}_i\rightarrow a^{(s+1)}_i)\right)=1$ by (\ref{dmdmdm}). Hence we assume $r-m'<i$. Adding (\ref{thm1lem-rel1}), we suppose that
\[
0= a^{(l^{(\delta_i+1)}_i)}_i=a^{(l^{(\delta_i+1)}_i+1)}_i=\cdots=a^{(\beta)}_i
<a^{(\beta+1)}_i=\cdots=a^{(\gamma)}_i=\ovl{r}<a^{(\gamma+1)}_i,
\]
for some numbers $\beta$, $\gamma$ $(l^{(\delta_i+1)}_i-1\leq\beta\leq\gamma\leq l^{(m-m')}_i)$. Using Definition \ref{labeldef}, if $l^{(\delta_i+1)}_i\leq\beta$ then
\begin{eqnarray*}
& &\prod^{l^{(m-m')}_i}_{s=l^{(\delta_i+1)}_i-1}\left(Q^{(s)}(a^{(s)}_i\rightarrow a^{(s+1)}_i)\right) \\
&=&\left(\frac{1}{Y_{m-l^{(\delta_i+1)}_i+1,r}}\right)\cdot
\left(\frac{2Y_{m-\beta,r}}{Y_{m-\beta,r-1}}\right)\cdot
\left(\prod^{\gamma-1}_{s=\beta+1}
\frac{Y^2_{m-s,r}}{Y_{m-s,r-1}}\right) \\
& &\qq \qq \qq \qq \qq \qq \cdot
\left(\frac{Y^2_{m-\gamma,r}}{Y_{m-\gamma,|a^{(\gamma+1)}_i|-1}}\right)\cdot
\left(\prod^{l^{(m-m')}_i}_{s=\gamma+1}
\frac{Y_{m-s,|a^{(s)}_i|}}{Y_{m-s,|a^{(s+1)}_i|-1}}\right)
\qq \qq \\
&=&2\left(\prod^{\beta}_{s=l^{(\delta_i+1)}_i} \frac{Y_{m-s,r}}{Y_{m-s+1,r}}\right)\cdot
\left(\prod^{\gamma}_{s=\beta+1}
\frac{Y^2_{m-s,r}}{Y_{m-s+1,r-1}} \right) \\
& &\qquad \qq \qq \qq \qq \qq \qq \cdot
\left(\prod^{l^{(m-m')}_i}_{s=\gamma+1}
\frac{Y_{m-s,|a^{(s)}_i|}}{Y_{m-s+1,|a^{(s)}_i|-1}}\right)\cdot\frac{1}{Y_{m'+i-r,d-i}} \\
&=&2\left(\prod^{m-m'}_{\zeta=\delta_i+1}B(m-l^{(\zeta)}_i,k^{(\zeta)}_i)\right)\cdot\frac{1}{Y_{m'+i-r,d-i}},
\end{eqnarray*}
where we used $l^{(m-m')}_i=m-i+r-m'$ and $|a^{(l^{(m-m')}_i+1)}_i|=d-i+1$ (Lemma \ref{cancellem}) in the second equality. Similarly, if $\beta=l^{(\delta_i+1)}_i-1$ then
\begin{equation*}
\prod^{l^{(m-m')}_i}_{s=l^{(\delta_i+1)}_i-1}\left(Q^{(s)}(a^{(s)}_i\rightarrow a^{(s+1)}_i)\right)=
\left(\prod^{m-m'}_{\zeta=\delta_i+1}B(m-l^{(\zeta)}_i,k^{(\zeta)}_i)\right)\cdot\frac{1}{Y_{m'+i-r,d-i}}.
\end{equation*}
Therefore, 
\begin{multline}\label{thm1lempr3-0}
\prod^{l^{(m-m')}_i}_{s=l^{(\delta_i+1)}_i-1} c^{(s)}(p)\prod^{d}_{i=1}\left(Q^{(s)}(a^{(s)}_i\rightarrow a^{(s+1)}_i)\right)=\\
2^{C[\{k^{(s)}_i\}]}\cdot\prod^{d}_{i=1}\left(\prod^{m-m'}_{\zeta=\delta_i+1}B(m-l^{(\zeta)}_i,k^{(\zeta)}_i)\right)\cdot\frac{1}{Y_{m'+i-r,d-i}},
\end{multline}
where $C[\{k^{(s)}_i\}]=C[\{k^{(s)}_i\}_{1\leq i\leq d,\ 1\leq s\leq m-m'}]$ is the number defined in (\ref{Cpdef}).

Finally, we consider the last range $(l^{(m-m')}_i+1\leq s\leq m-1)$ of the product. Using Lemma \ref{cancellem}, Lemma \ref{kpro2} and (\ref{ldef}), we obtain
\begin{equation}\label{thm1lempr4}
\prod^{m-1}_{s=l^{(m-m')}_i+1}\left(Q(a^{(s)}_i\rightarrow a^{(s+1)}_i)\right)=
\begin{cases}
\prod^{m-1}_{s=m-m'-i+r+1} \frac{Y_{m-s, d-i+1}}{Y_{m-s, d-i}} & {\rm if}\ r-m'<i, \\
1 & {\rm if}\ r-m'\geq i.
\end{cases}
\end{equation}
By (\ref{thm1lempr2}), (\ref{thm1lempr3-0}) and (\ref{thm1lempr4}), to prove (\ref{thm1lemclaim}), we need to show that
\begin{equation}\label{thm1lempr5}
\prod^d_{i=r-m'+1}\left(\frac{1}{Y_{m'+i-r, d-i}} \prod^{m-1}_{s=m-m'-i+r+1}\frac{Y_{m-s, d-i+1}}{Y_{m-s, d-i}}\right)=1. 
\end{equation}
We set 
\[A:=\prod^d_{i=r-m'+1}\left(\frac{1}{Y_{m'+i-r, d-i}}\right),\ {\rm and}\q B:=\prod^d_{i=r-m'+1}\left(\prod^{m-1}_{s=m-m'-i+r+1}\frac{Y_{m-s, d-i+1}}{Y_{m-s, d-i}}\right).\]

We obtain the followings:
\[
A=\prod^d_{i=r-m'+1}\left(\frac{1}{Y_{m'+i-r, d-i}}\right)=
\prod^{d-1}_{i=r-m'+1}\left(\frac{1}{Y_{m'+i-r, d-i}}\right)=
\prod^{m'+d-r-1}_{k=1}\left(\frac{1}{Y_{k, d-r+m'-k}}\right),
\]
and
\begin{eqnarray*}
B&=&\prod^d_{i=r-m'+1}\left(\prod^{m-1}_{s=m-m'-i+r+1}\frac{Y_{m-s, d-i+1}}{Y_{m-s, d-i}}\right)=
\prod^d_{i=r-m'+1}\left(\prod^{m'+i-r-1}_{s=1}\frac{Y_{s, d-i+1}}{Y_{s, d-i}}\right)\\
&=&\prod^{m'+d-r-1}_{s=1}\left(\frac{Y_{s,d-r+m'-s}}{Y_{s,d-r+m'-s-1}}\frac{Y_{s,d-r+m'-s-1}}{Y_{s,d-r+m'-s-2}}
\frac{Y_{s,d-r+m'-s-2}}{Y_{s,d-r+m'-s-3}}\cdots \frac{Y_{s,1}}{Y_{s,0}} \right)\\
&=&\prod^{m'+d-r-1}_{s=1} Y_{s,d-r+m'-s},
\end{eqnarray*}
where we used $Y_{s,0}=1$ (see Remark \ref{importantrem}). Thus we have $A\cdot B=1$,
which implies $(\ref{thm1lempr5})$. \qed

\vspace{3mm}

Let us prove the main theorem.

\nd
{\sl Proof of Theorem \ref{thm1}.}

As we have seen in Lemma \ref{thm1lem}, the monomial $Q(p)$ $(p\in X_d(m,m'))$ is described as (\ref{thm1lemclaim}) with $\{k^{(s)}_i\}_{1\leq i\leq d, 1\leq s\leq m-m'}$ which satisfies the conditions in (\ref{kpro1}) and Lemma \ref{kpro2} (ii) (iii), that is:
\[ 1\leq k^{(s)}_1\leq k^{(s)}_2\leq \cdots\leq k^{(s)}_d\leq\ovl{1},\q \ k^{(s)}_i=k^{(s)}_{i+1}\ {\rm if\ and\ only\ if}\ k^{(s)}_i=k^{(s)}_{i+1}=0,\]
\[
i\leq k^{(1)}_i\leq \cdots\leq k^{(m-m')}_i\leq m'+i \q (1\leq i\leq r-m'), \]
\[
i\leq k^{(1)}_i\leq \cdots\leq k^{(m-m')}_i\leq \ovl{d-i+1} \q (r-m'< i\leq d), \\
\]
\[
k^{(s)}_i<k^{(s+1)}_{i+1}. \\
\]
Thus, $\{k^{(s)}_i\}$ satisfies the conditions $(*)$ in Theorem \ref{thm1}.

Let $\{K^{(s)}_i\}_{1\leq i\leq d, 1\leq s\leq m-m'}$ be the set of numbers which satisfies the conditions $(*)$ in Theorem \ref{thm1}:
\begin{equation}\label{cond1}
1\leq K^{(s)}_1\leq K^{(s)}_2\leq \cdots\leq K^{(s)}_d\leq\ovl{1},\ K^{(s)}_i=K^{(s)}_{i+1}\ {\rm if\ and\ only\ if}\ K^{(s)}_i=K^{(s)}_{i+1}=0,
\end{equation}
\begin{equation}\label{cond2}
i\leq K^{(1)}_i\leq \cdots\leq K^{(m-m')}_i\leq m'+i\q (1\leq i\leq r-m'),
\end{equation}
\begin{equation}\label{cond3}
i\leq K^{(1)}_i\leq \cdots\leq K^{(m-m')}_i\leq \ovl{1}\q (r-m'< i\leq d),
\end{equation}
and
\begin{equation}\label{cond4}
K^{(s)}_i<K^{(s+1)}_{i+1}.
\end{equation}

We need to show that there exists a path $p\in X_d(m,m')$ such that
\begin{equation}\label{finclaim}
Q(p)= 2^{C[\{K^{(s)}_i\}]}\prod^{d}_{i=1} B(m-L^{(1)}_i,K^{(1)}_i)\cdots B(m-L^{(m-m')}_i,K^{(m-m')}_i),
\end{equation}
where $C[\{K^{(s)}_i\}]$ is defined as in (\ref{Cpdef}) and
\[ L^{(s)}_{i}:=
\begin{cases}
K^{(s)}_{i}+s-i-1 & {\rm if}\  K^{(s)}_{i}\in\{j|1\leq j\leq r\}, \\
s-i+r & {\rm if}\ K^{(s)}_{i}\in\{\ovl{j}|1\leq j\leq r\}\cup\{0\},
\end{cases}
\]
for $1\leq s\leq m-m'$ and $1\leq i\leq d$. Since we supposed $K^{(s)}_i\leq K^{(s)}_{i+1}$ and $K^{(s)}_i=K^{(s)}_{i+1}$ if and only if $K^{(s)}_i=K^{(s)}_{i+1}=0$, we can easily verify \begin{equation}\label{lequ1}
L^{(s)}_i\leq L^{(s)}_{i+1}\q {\rm if}\ K^{(s)}_i\in\{j|1\leq j\leq r\},
\end{equation}
and 
\begin{equation}\label{lequ2}
L^{(s)}_i=L^{(s)}_{i+1}+1\q {\rm if}\ K^{(s)}_i\in\{\ovl{j}|1\leq j\leq r\}\cup\{0\}. 
\end{equation}

We claim that $0\leq L^{(s)}_i\leq m-1$. By the conditions (\ref{cond2}) and (\ref{cond3}), we get $i\leq K^{(s)}_i$. So it is clear that $0\leq L^{(s)}_i$. For $i$ and $s$ $(1\leq i\leq r-m',\ 1\leq s\leq m-m')$, it follows from the condition $(\ref{cond2})$ that $L^{(s)}_{i}=K^{(s)}_{i}+s-i-1\leq m'+i+s-i-1=m'+s-1\leq m-1$. For $i$ $(r-m'<i)$, we get $L^{(s)}_i\leq r-i+s<m'+s\leq m$. Therefore, we have $0\leq L^{(s)}_i\leq m-1$ for all $1\leq i\leq d$ and $1\leq s\leq m-m'$.

Note that if $K^{(s)}_i\in\{\ovl{j}|1\leq j\leq r\}\cup\{0\}$, then $0\leq K^{(s)}_i\leq K^{(s)}_{i+1}\leq \ovl{1}$ and hence
\begin{equation}\label{lequ3}
L^{(s+1)}_i=L^{(s)}_i+1.
\end{equation}
We define a path $p={\rm vt}(m;a^{(0)}_1,\cdots,a^{(0)}_d)\rightarrow\cdots \rightarrow {\rm vt}(0;a^{(m)}_1,\cdots,a^{(m)}_d)\in X_d(m,m')$ as follows: For $i$ $(1\leq i\leq r-m')$, we define the $i$-sequence (Definition \ref{iseq}) of $p$ as
\begin{flushleft}
$a^{(0)}_i=i,\ a^{(1)}_i=i+1,\ a^{(2)}_i=i+2,\cdots,a^{(L^{(1)}_i)}_i=i+L^{(1)}_i$,
 
$ a^{(L^{(1)}_i+1)}_i=i+L^{(1)}_i,\ a^{(L^{(1)}_i+2)}_i=i+L^{(1)}_i+1,\cdots, a^{(L^{(2)}_i)}_i=i+L^{(2)}_i-1, $

$ a^{(L^{(2)}_i+1)}_i=i+L^{(2)}_i-1,\ a^{(L^{(2)}_i+2)}_i=i+L^{(2)}_i,\cdots, a^{(L^{(3)}_i)}_i=i+L^{(3)}_i-2, $
\end{flushleft}
\begin{equation}\label{jlist2} \vdots \end{equation}
\begin{flushleft}
$a^{(L^{(m-m'-1)}_i+1)}_i=i+L^{(m-m'-1)}_i-m+m'+2,\ 
\cdots, a^{(L^{(m-m')}_i)}_i=i+L^{(m-m')}_i-m+m'+1,$

$a^{(L^{(m-m')}_i+1)}_i=i+L^{(m-m')}_i-m+m'+1,\ a^{(L^{(m-m')}_i+2)}_i=i+L^{(m-m')}_i-m+m'+2,$

$a^{(L^{(m-m')}_i+3)}_i=i+L^{(m-m')}_i-m+m'+3,\cdots, a^{(m)}_i=m'+i$.
\end{flushleft}

\vspace{2mm}

For $i$ $(r-m'+1\leq i\leq d)$, we define the $i$-sequence of $p$ as follows: We suppose that $1\leq K^{(1)}_i\leq \cdots\leq K^{(\delta_i)}_i\leq r<0\leq K^{(\delta_i+1)}_i\leq \cdots \leq K^{(m-m')}_i\leq\ovl{1}$ for some $\delta_i$ $(1\leq\delta_i\leq m-m')$. Then
\begin{flushleft}
$a^{(0)}_i=i,\ a^{(1)}_i=i+1,\ a^{(2)}_i=i+2,\cdots,a^{(L^{(1)}_i)}_i=i+L^{(1)}_i$,
 
$ a^{(L^{(1)}_i+1)}_i=i+L^{(1)}_i,\ a^{(L^{(1)}_i+2)}_i=i+L^{(1)}_i+1,\cdots, a^{(L^{(2)}_i)}_i=i+L^{(2)}_i-1, $

$ a^{(L^{(2)}_i+1)}_i=i+L^{(2)}_i-1,\ a^{(L^{(2)}_i+2)}_i=i+L^{(2)}_i,\cdots, a^{(L^{(3)}_i)}_i=i+L^{(3)}_i-2, $
\end{flushleft}
\begin{equation}\label{jlist3} \vdots \end{equation}
\begin{flushleft}
$a^{(L^{(\delta_i-1)}_i+1)}_i=i+L^{(\delta_i-1)}_i-\delta_i+2,\ 
\cdots, a^{(L^{(\delta_i)}_i)}_i=i+L^{(\delta_i)}_i-\delta_i+1,$

$a^{(L^{(\delta_i)}_i+1)}_i=i+L^{(\delta_i)}_i-\delta_i+1,\ a^{(L^{(\delta_i)}_i+2)}_i=i+L^{(\delta_i)}_i-\delta_i+2,$

$a^{(L^{(\delta_i)}_i+3)}_i=i+L^{(\delta_i)}_i-\delta_i+3,\cdots,\ a^{(L^{(\delta_i+1)}_i-1)}_i=r$,

$a^{(L^{(\delta_i+1)}_i)}_i=K^{(\delta_i+1)}_i,\ a^{(L^{(\delta_i+2)}_i)}_i=K^{(\delta_i+2)}_i,\cdots,\ a^{(L^{(m-m')}_i)}_i=K^{(m-m')}_i$, 

$a^{(L^{(m-m')}_i+1)}_i=a^{(L^{(m-m')}_i+2)}_i=\cdots=a^{(m)}_i=\ovl{d-i+1}$.
\end{flushleft}

\vspace{2mm}

Above lists say $a^{(L^{(s)}_i)}_i=K^{(s)}_i$ $(1\leq s\leq m-m')$. Clearly, the path $p$ satisfies the conditions in Definition \ref{pathdef} (iii) and (iv). For $s$ $(1\leq s\leq L^{(\delta_{i}+1)}_i-1)$, it follows from (\ref{lequ1}) and the lists (\ref{jlist2}), (\ref{jlist3}) that $a^{(s)}_i<a^{(s)}_{i+1}$. For $s$ $(\delta_{i}+1\leq s\leq m-m')$, we obtain $a^{(L^{(s)}_i)}_i<a^{(L^{(s)}_{i})}_{i+1}$ since  $a^{(L^{(s)}_i)}_i=K^{(s)}_i<K^{(s+1)}_{i+1}=a^{(L^{(s+1)}_{i+1})}_{i+1}=a^{(L^{(s)}_{i+1}+1)}_{i+1}=
a^{(L^{(s)}_{i})}_{i+1}$ by (\ref{cond4}), (\ref{lequ2}) and (\ref{lequ3}). For $s$ $(L^{(m-m')}_i+1\leq s\leq m)$, we obtain $a^{(s)}_i=\ovl{d-i+1}$, and $a^{(s)}_{i+1}=\ovl{d-i}$ since $L^{(m-m')}_{i+1}=L^{(m-m')}_i-1<L^{(m-m')}_i\leq s$, which means $a^{(s)}_i<a^{(s)}_{i+1}$. Therefore, $a^{(s)}_i<a^{(s)}_{i+1}$ for all $1\leq i\leq d-1$ and $1\leq s\leq m-m'$, which means the path $p$ satisfies Definition \ref{pathdef} (ii).

Finally, for $a^{(s)}_i\in \{\ovl{j}|1\leq j\leq r\}\cup\{0\}$, we need to verify $a^{(s)}_i\leq a^{(s-1)}_{i+1}$ and $a^{(s)}_i=a^{(s-1)}_{i+1}$ if and only if $a^{(s)}_i=a^{(s-1)}_{i+1}=0$. The list (\ref{jlist3}) of $i$-sequence of $p$ shows that either $s=L^{(\zeta)}_i$ for some $\zeta$ $(\delta_i+1\leq\zeta\leq m-m')$ or $L^{(m-m')}_i<s$ holds. If $s=L^{(\zeta)}_i$, using (\ref{cond1}) and (\ref{lequ2}), we see that $a^{(s)}_i=a^{(L^{(\zeta)}_i)}_i=K^{(\zeta)}_i\leq K^{(\zeta)}_{i+1}=
a^{(L^{(\zeta)}_{i+1})}_{i+1}=a^{(L^{(\zeta)}_{i}-1)}_{i+1}=a^{(s-1)}_{i+1}$, and $a^{(s)}_i=a^{(s-1)}_{i+1}$ if and only if $a^{(s)}_i=a^{(s-1)}_{i+1}=0$. In the case $L^{(m-m')}_i<s$, we obtain $a^{(s)}_i=\ovl{d-i+1}<\ovl{d-i}=a^{(s-1)}_{i+1}$ since $L^{(m-m')}_{i+1}=L^{(m-m')}_{i}-1<L^{(m-m')}_{i}\leq s-1$. Therefore, we have $a^{(s)}_i<a^{(s-1)}_{i+1}$ for $a^{(s)}_i\in \{\ovl{j}|1\leq j\leq r\}\cup\{0\}$, which means the path $p$ satisfies Definition \ref{pathdef} (v). Hence $p$ is well-defined, and (\ref{finclaim}) follows from Lemma \ref{thm1lem}, and Theorem \ref{thm1} follows from Proposition \ref{pathlem}. \qed

\subsection{The proof of Theorem \ref{thm2}}

In the case $i_k=r$, explicit forms and path descriptions of the minors $\Delta^L(k;\textbf{i})$ are different from the ones in Theorem \ref{thm1} and Proposition \ref{pathlem}. First, we shall prove the path description of $\Delta^L(k;\textbf{i})$ for $i_k=r$.

\begin{defn}\label{vspesp}
Let us define the directed graph $(V^{sp},E^{sp})$ as follows:
We define the set $V^{sp}=V^{sp}(m)$ of vertices as 
\begin{multline*}
V^{sp}(m):=\{{\rm vt}(s;k^{(s)}_1,k^{(s)}_2,\cdots,k^{(s)}_t)|0\leq s\leq m,\\
 k^{(s)}_i\in \{1,2,\cdots,r\},\ t\leq s,\ k^{(s)}_1<k^{(s)}_2<\cdots<k^{(s)}_t \}. 
\end{multline*}
And we define the set $E^{sp}=E^{sp}(m)$ of directed edges as 
\begin{multline*} E^{sp}(m):=\{{\rm vt}(s;k^{(s)}_1,\cdots,k^{(s)}_t)\rightarrow
{\rm vt}(s-1;k^{(s+1)}_1,\cdots,k^{(s+1)}_{t'})\\
|\ 0\leq s\leq m-1,\ {\rm vt}(s;k^{(s)}_1,\cdots,k^{(s)}_t),\ {\rm vt}(s-1;k^{(s-1)}_1,\cdots,k^{(s-1)}_{t'})\in V^{sp}(m)\}.
\end{multline*}
\end{defn}

Now, let us define the set of directed paths from ${\rm vt}(m;\ )$ to ${\rm vt}(0;1,2,\cdots,m')$ in $(V^{sp},E^{sp})$.
\begin{defn}\label{pathdef-sp}
Let $X^r(m,m')$ be the set of directed paths $p$ from ${\rm vt}(m;\ )$ to ${\rm vt}(0;1,2,\cdots,m')$ in $(V^{sp},E^{sp})$ which satisfy the following conditions: For $s\in\mathbb{Z}$ $(1\leq s\leq m)$, if ${\rm vt}(s;k^{(s)}_1,\cdots,k^{(s)}_t)\rightarrow{\rm vt}(s-1;k^{(s-1)}_1,\cdots,k^{(s-1)}_{t'})$ is an edge included in $p$, then
\begin{enumerate}
\item $t'=t$ or $t'=t+1$, if $k^{(s)}_t=r$ then $t'=t$,
\item $k^{(s-1)}_i\leq k^{(s)}_i<k^{(s-1)}_{i+1}$,
\item $k^{(0)}_i=i$ $(1\leq i\leq m')$.
\end{enumerate}
\end{defn}

\begin{defn}\label{connected-def-sp}
If two vertices ${\rm vt}(s;k^{(s)}_1,\cdots,k^{(s)}_t)$ and ${\rm vt}(s-1;k^{(s-1)}_1,\cdots,k^{(s-1)}_{t'})$ satisfy the conditions (i) and (ii) in Definition \ref{pathdef-sp}, then we say that these vertices are {\it connected}. 
\end{defn}

Define a Laurent monomial associated with each edge of paths in $X^r(m,m')$.

\begin{defn}\label{labeldef-sp}
Let $p\in X^r(m,m')$ be a path.
\item For each $s$ $(0\leq s\leq m)$, we define the {\it label of the edge} ${\rm vt}(s;i_1,i_2,\cdots,i_t)\rightarrow{\rm vt}(s-1;j_1,j_2,\cdots,j_{t'})$ as the Laurent monomial which is given as follows and write it as $Q^{(s)}(p)$ or $Q({\rm vt}(s;i_1,i_2,\cdots,i_t)\rightarrow{\rm vt}(s-1;j_1,j_2,\cdots,j_{t'}))$: 
\[ Q^{(s)}(p):=
\begin{cases}
\frac{Y_{m-s,i_1}}{Y_{m-s,j_1-1}}\frac{Y_{m-s,i_2}}{Y_{m-s,j_2-1}}\cdots \frac{Y_{m-s,i_t}}{Y_{m-s,j_t-1}}\cdot \frac{1}{Y_{m-s,r}} & {\rm if}\ t'=t, \\
\frac{Y_{m-s,i_1}}{Y_{m-s,j_1-1}}\frac{Y_{m-s,i_2}}{Y_{m-s,j_2-1}}\cdots \frac{Y_{m-s,i_t}}{Y_{m-s,j_t-1}}\frac{1}{Y_{m-s,j_{t+1}-1}} & {\rm if}\ t'=t+1. 
\end{cases}
.\]
And we define the {\it label} $Q(p)$ {\it of the path} $p$ as the total product:
\[Q(p):=\prod_{s=0}^{m-1}Q^{(s)}(p). 
\]
\end{defn}

\begin{prop}\label{pathlem-spin}In the setting of Theorem \ref{thm2},
we have the following:
\begin{equation}\label{pathlem-spin-e}
 \Delta^L(k;{\rm \bf{i}})(\textbf{Y})=\sum_{p\in X^r(m,m')} Q(p). 
\end{equation}
\end{prop}

{\sl Proof.}

Using the bilinear form (\ref{minor-bilin}), we obtain
\[ \Delta^L(k;\textbf{i})(\textbf{Y})=
 \Del_{u_{\leq k}\Lm_r,\Lm_r}(x^L_{\textbf{i}}(\textbf{Y}))
=\lan x^L_{\textbf{i}}(\textbf{Y})\cdot v_{\Lm_r}\, ,\, \ovl{u_{\leq k}}\cdot v_{\Lm_r}\ran,
\]
where $v_{\Lm_r}$ is the highest vector in the spin representation $V(\Lm_r)$. As we have seen in Sect.\ref{SectFundB}, we can describe the weight vectors of $V(\Lm_r)$ as the form $(\ep_1,\cdots,\ep_r)$ $(\ep_i\in\{+,-\})$. In particular, the highest weight vector $v_{\Lm_r}$ is equal to $(+,+,\cdots,+)$. Using the definitions (\ref{Bsp-f1}) and (\ref{Bsp-f2}), we have
\begin{eqnarray*}
\ovl{u_{\leq k}}\cdot v_{\Lm_r}&=&(\ovl{s_1}\ \ovl{s_2}\cdots\ovl{s_r})^{m'}(+,+,+,\cdots,+) \\
&=&(\ovl{s_1}\ \ovl{s_2}\cdots\ovl{s_r})^{m'-1}(-,+,+,\cdots,+) \\
&\vdots&\\
&=&(\underbrace{-,-,\cdots,-}_{m'},+,\cdots,+).
\end{eqnarray*}
Now, we define the following notation: For $(\ep_1,\cdots,\ep_r)\in{\mathbf B}_{{\rm sp}}^{(r)}$, we suppose that
\[ \{1\leq i\leq r|\ \ep_i=- \}=\{k_1,\cdots,k_t\},\ \ \ (k_1<\cdots<k_t). \]
Then we define $[k_1,\cdots,k_t]:=(\ep_1,\cdots,\ep_r)$. Acting $x_{-1}(Y_{s,1})\cdots x_{-r}(Y_{s,r})$ on $[k_1,\cdots,k_t]$ $(1\leq s\leq m)$, it becomes a linear combination of $\{[j_1,\cdots,j_t]|\ k_{l-1}<j_l\leq k_l\ (1\leq l\leq t)\}\cup\{[j_1,\cdots,j_t,j_{t+1}]|k_{l-1}<j_l\leq k_l\ (1\leq l\leq t+1)\}$. It follows from the definitions (\ref{Bsp-f0}) and (\ref{Bsp-f1}) of actions that
\begin{multline*}
x_{-k_{l-1}}(Y_{s,k_{l-1}})\cdots x_{-(k_l-1)}(Y_{s,k_l-1})[\cdots,k_l,\cdots]=\cdots+\left(\frac{1}{Y_{s,k_l-1}}\right)^{\frac{1}{2}}\left(\frac{Y_{s,k_l-1}}{Y_{s,k_l-2}}\right)^{\frac{1}{2}}\cdots \\
\left(\frac{Y_{s,j_l}}{Y_{s,j_l-1}}\right)^{\frac{1}{2}} \left(\frac{Y_{s,j_l-2}}{Y_{s,j_l-1}}\right)^{\frac{1}{2}}\left(\frac{Y_{s,j_l-3}}{Y_{s,j_l-2}}\right)^{\frac{1}{2}}\cdots \left(\frac{Y_{s,k_{l-1}+1}}{Y_{s,k_{l-1}+2}}\right)^{\frac{1}{2}}\frac{Y_{s,k_{l-1}}}{(Y_{s,k_{l-1}+1})^{\frac{1}{2}}}[\cdots,j_l,\cdots]+\cdots \\
=\cdots+\frac{Y_{s,k_{l-1}}}{Y_{s,j_l-1}}[\cdots,j_l,\cdots]+\cdots\qq \qq \qq \qq \qq \qq.
\end{multline*}
Thus, the coefficient of $[j_1,\cdots,j_t]$ in $x_{-1}(Y_{s,1})\cdots x_{-r}(Y_{s,r})[k_1,\cdots,k_t]$ is
\[ \frac{Y_{s,k_1}}{Y_{s,j_1-1}}\frac{Y_{s,k_2}}{Y_{s,j_2-1}}\cdots \frac{Y_{s,k_t}}{Y_{s,j_t-1}}\cdot \frac{1}{Y_{s,r}}, \]
and the one of $[j_1,\cdots,j_t,j_{t+1}]$ is
\[ \frac{Y_{s,k_1}}{Y_{s,j_1-1}}\frac{Y_{s,k_2}}{Y_{s,j_2-1}}\cdots \frac{Y_{s,k_t}}{Y_{s,j_t-1}}\frac{1}{Y_{s,j_{t+1}-1}}, \]
which coincides with $Q({\rm vt}(m-s;k_1,\cdots,k_t)\rightarrow{\rm vt}(m-s+1;j_1,\cdots,j_{t'}))$ in Definition \ref{labeldef-sp}. Note that if $k_t=r$ then $x_{-1}(Y_{s,1})\cdots x_{-r}(Y_{s,r})[k_1,\cdots,k_t]$ is a linear combination of $\{[j_1,\cdots,j_t]|\ k_{l-1}<j_l\leq k_l\ (1\leq l\leq t)\}$. In our notation, we can write $[\ ]=(+,+,\cdots,+)$ and $[1,2,\cdots,m']=(\underbrace{-,-,\cdots,-}_{m'},+,\cdots,+)$. Therefore, we get (\ref{pathlem-spin-e}). \qed

\vspace{2mm}

{\sl Proof of Theorem \ref{thm2}.}

For a path $p\in X^r(m,m')$, we write it explicitly as
\begin{multline*}{\rm vt}(m;\ )\rightarrow{\rm vt}(m-1;\ )\rightarrow \cdots 
\rightarrow{\rm vt}(t_{1};\ )\rightarrow
{\rm vt}(t_{1}-1; k^{(t_{1}-1)}_1)\rightarrow\cdots
\rightarrow \\
{\rm vt}(t_{2}; k^{(t_{2})}_1)\rightarrow
{\rm vt}(t_{2}-1; k^{(t_{2}-1)}_1,k^{(t_{2}-1)}_2)\rightarrow\cdots\\
\rightarrow
{\rm vt}(t_{m'}; k^{(t_{m'})}_1,\cdots,k^{(t_{m'})}_{m'-1})
\rightarrow
{\rm vt}(t_{m'}-1; k^{(t_{m'}-1)}_1,\cdots,k^{(t_{m'}-1)}_{m'-1},k^{(t_{m'}-1)}_{m'})
\\
\rightarrow \cdots \rightarrow
{\rm vt}(1; k^{(1)}_1,k^{(1)}_2,\cdots, k^{(1)}_{m'})\rightarrow
{\rm vt}(0; 1,2,\cdots,m')
\end{multline*}
with some integers $\{t_i\}^{m'}_{i=1}$ and $\{k^{(s)}_i\}^{t_i-1}_{s=1}$ such that $1\leq t_{m'}<t_{m'-1}<\cdots<t_1\leq m,\ i\leq k^{(1)}_i\leq \cdots \leq k^{(t_i-1)}_i\leq r$, $k^{(s-1)}_i\leq k^{(s)}_i<k^{(s-1)}_{i+1}$ $(1\leq s\leq t_i-1)$. The previous proposition implies
\begin{multline*} Q(p)=\left(\prod^{m'}_{i=1}\frac{1}{Y_{t_i,k^{(t_i-1)}_i-1}}\frac{Y_{t_i-1,k^{(t_i-1)}_i}}{Y_{t_i-1,k^{(t_i-2)}_i-1}}\frac{Y_{t_i-2,k^{(t_i-2)}_i}}{Y_{t_i-2,k^{(t_i-3)}_i-1}}\cdots \frac{Y_{1,k^{(1)}_i}}{Y_{1,k^{(0)}_i-1}} \right) \\
\qq \qq  \cdot\left(\frac{1}{Y_{m,r}Y_{m-1,r}\cdots Y_{t_1+1,r}Y_{t_1-1,r}\cdots Y_{t_2+1,r}Y_{t_2-1,r}\cdots Y_{t_{m'}+1,r}Y_{t_{m'}-1,r}\cdots Y_{1,r}} \right).
\end{multline*}
It can be easily seen that $k^{(m'-i)}_i=k^{(m'-i-1)}_i=\cdots=k^{(0)}_i=i$ by Definition \ref{pathdef-sp} (ii) and (iii). Hence, setting $t_0:=m+1$ and $t_{m'+1}:=0$, we have
\begin{eqnarray}\label{thm2-pr1}
Q(p)&=&\prod^{m'+1}_{i=1}\left(\frac{1}{Y_{t_i,k^{(t_i-1)}_i-1}}\frac{Y_{t_i-1,k^{(t_i-1)}_i}}{Y_{t_i-1,k^{(t_i-2)}_i-1}}\cdots \frac{Y_{m'-i+2,k^{(m'-i+2)}_i}}{Y_{m'-i+2,k^{(m'-i+1)}_i-1}} \frac{Y_{m'-i+1,k^{(m'-i+1)}_i}}{Y_{m'-i+1,i-1}} \right) \nonumber
\\ \nonumber
& & \cdot\left(\frac{Y_{m'-i,i}}{Y_{m'-i,i-1}}\frac{Y_{m'-i-1,i}}{Y_{m'-i-1,i-1}}\cdots \frac{Y_{1,i}}{Y_{1,i-1}}\right)
\cdot\left(\frac{1}{Y_{t_{i-1}-1,r}}\frac{1}{Y_{t_{i-1}-2,r}}\cdots \frac{1}{Y_{t_{i}+1,r}} \right)\\ \nonumber
&=&\prod^{m'+1}_{i=1}\left(\frac{Y_{t_i-1,k^{(t_i-1)}_i}}{Y_{t_i,k^{(t_i-1)}_i-1}}\frac{Y_{t_i-2,k^{(t_i-2)}_i}}{Y_{t_i-1,k^{(t_i-2)}_i-1}}\cdots\frac{Y_{m'-i+1,k^{(m'-i+1)}_i}}{Y_{m'-i+2,k^{(m'-i+1)}_i-1}} \right)\\ \nonumber
& &\qq \qq \qq \qq \qq \qq \qq
\cdot\left(\frac{1}{Y_{t_{i-1}-1,r}}\frac{1}{Y_{t_{i-1}-2,r}}\cdots \frac{1}{Y_{t_{i}+1,r}} \right)\\ \nonumber
&=&\prod^{m'+1}_{i=1} B(t_i-1,\ovl{k^{(t_i-1)}_i})B(t_i-2,\ovl{k^{(t_i-2)}_i})\cdots B(m'-i+1,\ovl{k^{(m'-i+1)}_i})\\ 
& & \cdot
B(t_{i-1}-1,r+1)B(t_{i-1}-2,r+1)\cdots B(t_{i}+1,r+1),
\end{eqnarray}
where in the second equality, we used the following: 
\[\prod^{m'+1}_{i=1}\left(\frac{1}{Y_{m'-i+1,i-1}}\cdot\frac{Y_{m'-i,i}}{Y_{m'-i,i-1}}\frac{Y_{m'-i-1,i}}{Y_{m'-i-1,i-1}}\cdots \frac{Y_{1,i}}{Y_{1,i-1}}\right)=1. \]
Therefore, each term of $\Delta^L(k;\textbf{i})(\textbf{Y})$ is described as (\ref{thm2-pr1}). 

Conversely, if integers $\{T_i\}^{m'}_{i=1}$ and $\{K^{(s)}_i\}^{T_i-1}_{s=m'-i+1}$ such that $1\leq T_{m'}<T_{m'-1}<\cdots<T_1\leq m,\ i\leq K^{(m'-i+1)}_i\leq \cdots \leq K^{(T_i-1)}_i\leq r$, $K^{(s-1)}_i\leq K^{(s)}_i<K^{(s-1)}_{i+1}$, $K^{(0)}_i=i$ $(1\leq i\leq m')$, $T_0:=m+1$ and $T_{m'+1}=0$ are given, we can define a path $P\in X^r(m,m')$ as
\begin{multline*}{\rm vt}(m;\ )\rightarrow{\rm vt}(m-1;\ )\rightarrow \cdots \rightarrow
{\rm vt}(T_{1}; )\rightarrow
{\rm vt}(T_{1}-1; K^{(T_{1}-1)}_1)\rightarrow
\cdots\rightarrow \\
{\rm vt}(T_{2}; K^{(T_{2})}_1)\rightarrow
{\rm vt}(T_{2}-1; K^{(T_{2}-1)}_1,K^{(T_{2}-1)}_2)\rightarrow\cdots\\
\rightarrow
{\rm vt}(T_{3}; K^{(T_{3})}_1,K^{(T_{3})}_2)
\rightarrow
{\rm vt}(T_{3}-1; K^{(T_{3}-1)}_1,K^{(T_{3}-1)}_2,K^{(T_{3}-1)}_3)
\rightarrow \cdots \rightarrow \\
{\rm vt}(1; 1,2,\cdots,m'-1, K^{(1)}_{m'})\rightarrow
{\rm vt}(0; 1,2,\cdots,m').
\end{multline*}

In the same way as (\ref{thm2-pr1}), we can verify
\begin{eqnarray*}
Q(P)&=&\prod^{m'+1}_{i=1} B(T_i-1,\ovl{K^{(T_i-1)}_i})B(T_i-2,\ovl{K^{(T_i-2)}_i})\cdots B(m'-i+1,\ovl{K^{(m'-i+1)}_i})\\ 
& &\qq \qq \cdot
B(T_{i-1}-1,r+1)B(T_{i-1}-2,r+1)\cdots B(T_{i}+1,r+1).
\end{eqnarray*}
Thus, we obtain Theorem \ref{thm2}. \qed

\end{document}